\documentclass[11pt]{article}
\usepackage{enumerate}
\usepackage{booktabs}
\usepackage{rotating}
\usepackage{amssymb,a4wide,latexsym,makeidx,epsfig,fleqn}
\usepackage{amsthm}
\usepackage{amsmath}
\usepackage{enumerate}

\newtheorem{Exa}{Example}[section]

\begin{document}
\textwidth 180mm \textheight 225mm
\title{A modified generalized shift-splitting preconditioner for nonsymmetric saddle point problems
\thanks{Supported by the National Natural Science Foundation of China~(No. 11171273) and Innovation Foundation for Doctor Dissertation of Northwestern
Polytechnical University (No. CX201628).}}
\author{{Zhengge Huang, Ligong Wang\footnote{Corresponding author.} , Zhong Xu and Jingjing Cui}\\
{\small Department of Applied Mathematics, School of Science, Northwestern
Polytechnical University,}\\ {\small  Xi'an, Shaanxi 710072,
People's Republic
of China.}\\
{\small E-mails: ZhenggeHuang@mail.nwpu.edu.cn; lgwang@nwpu.edu.cn(or lgwangmath@163.com);}\\
{\small zhongxu@nwpu.edu.cn; JingjingCui@mail.nwpu.edu.cn}\\
 }
\date{}
\maketitle
\begin{center}
\begin{minipage}{130mm}
\vskip 0.3cm
\begin{center}
{\small {\bf Abstract}}
\end{center}
{\small For the nonsymmetric saddle point problems with nonsymmetric positive definite (1,1) parts, the modified generalized shift-splitting (MGSSP) preconditioner as well as the MGSSP iteration method are derived in this paper, which generalize the MSSP preconditioner and the MSSP iteration method newly developed by Huang and Su (J. Comput. Appl. Math. 2017), respectively. The convergent and semi-convergent analysis of the MGSSP iteration method are presented, and we prove that this method is unconditionally convergent and semi-convergent. In addition, some spectral properties of the preconditioned matrix are carefully analyzed. Numerical results demonstrate the robustness and effectiveness of the MGSSP preconditioner and the MGSSP iteration method, and also illustrate that the MGSSP iteration method outperforms the GSS and GMSS iteration methods, and the MGSSP preconditioner is superior to the shift-splitting (SS), generalized SS (GSS), modified SS (MSS) and generalized MSS (GMSS) preconditioners for the GMRES method for solving the nonsymmetric saddle point problems.

 }
\end{minipage}
\end{center}

\section{Introduction }
\label{sec:ch6-introduction}

In a wide variety of scientific and engineering applications,
such as mixed finite element approximation of elliptic partial differential
equations, the image reconstruction and registration, computational fluid dynamics, weighted least-squares problems, networks computer graphics, constrained optimization and son on \cite{2,3,53}, we need to solve the following nonsymmetric saddle point problems of the form
\begin{eqnarray}
\mathcal{A}u=\left(
    \begin{array}{cc}
      A & B \\
      -B^{T} & 0 \\
    \end{array}
  \right)\left(
           \begin{array}{c}
             x \\
             y \\
           \end{array}
         \right)=\left(
                   \begin{array}{c}
                     f \\
                     -g \\
                   \end{array}
                 \right)\equiv{b},
\end{eqnarray}
where $A\in{\mathbb{R}^{m\times{m}}}$ is nonsymmetric positive definite, $B\in{\mathbb{R}^{m\times{n}}}$ is a rectangular matrix, $p\in\mathbb{R}^{m}$ and $q\in\mathbb{R}^{n}$ are given vectors, with $n\leq{m}$. Here, $B^{T}$ denotes the transpose of $B$. The system of linear equations (1) is also termed
as a Karush-Kuhn-Tucker (KKT) system, or an augmented system \cite{56,57}. For a wider class of saddle point problems, the readers can refer to \cite{1}.

Since the matrices $A$ and $B$ are large and sparse in general, the iteration methods are often much more suitable for solving it than
direct methods. When $B$ is of full column rank, a large variety of effective iterative methods based on matrix splitting as well
as their numerical properties have been investigated in the literature. For example, Golub et al. \cite{6} developed the SOR-like method, and in the sequel, Bai et al. \cite{4,5} extended the SOR-like method to the generalized SOR (GSOR) method and the parameterized inexact Uzawa method, respectively. For SOR-like methods established recently, see \cite{7,58}. Based on the Uzawa method presented by Bramble et al. and Elman and Golub in \cite{8,9}, Bai et al. \cite{4,5}, Dai et al. \cite{10} and Ma and Zheng \cite{50} employed the Uzawa-type methods and so forth in recent years. Besides, Bai et al. put forward the well-known Hermitian and skew-Hermitian splitting (HSS) methods \cite{17} and its variants \cite{16,18,20,21}. On the basis of the shift-splitting (SS) of a non-Hermitian matrix \cite{40}, Cao et al. \cite{35} derived the SS iteration method as well as the SS preconditioner for the nonsingular saddle point problems, and Chen and Ma \cite{37} and Cao et al. \cite{38} generalized the SS iteration method and obtained the generalized SS (GSS) iteration method. To increase the convergence rate of the GSS iteration method, Huang and Su \cite{51} newly developed the modified shift-splitting (MSSP) iteration method.

If $B$ in (1) is rank deficient, then the coefficient matrix $\mathcal{A}$ in (1) is singular, and we call (1) the singular saddle point problem. Some iteration methods and preconditioning techniques for solving singular saddle point problems have been proposed in the recent literature, see, e.g., \cite{25,27,28,29}. Zheng et al. \cite{26} proposed some sufficient conditions for the semi-convergence of the GSOR method and determined the optimal iteration parameters. Bai \cite{31} derived some necessary and sufficient conditions to assure the semi-convergence of the HSS method. Chen et al. \cite{30} and Cao et al. \cite{44} investigated the generalized shift-splitting iteration method for singular saddle point problems. Very recently, Dou et al. \cite{42} introduced the modifying the parameterized inexact Uzawa (PIU) for singular saddle point problems, and Zheng and Lu \cite{41} proved the semi-convergence of the upper and lower triangular (ULT) splitting iteration method for singular saddle point problems.

Recently, based on the preconditioner \cite{40} studied for a class of non-Hermtian positive definite linear systems, Cao et al. \cite{35} presented a shift-splitting (SS) preconditioner
\begin{eqnarray}
\mathcal{P}_{SS}=\frac{1}{2}\left(
                         \begin{array}{cc}
                           \alpha I+A & B \\
                           -B^{T} & \alpha I \\
                         \end{array}
                       \right)\nonumber
\end{eqnarray}
for the saddle point problem (1), where $\alpha$ is a positive constant and $I$ is the identity matrix. The authors also proved the corresponding SS iteration method is unconditional convergent.

On the basis of the shift-splitting (SS) preconditioner \cite{35}, Chen and Ma \cite{37} and Cao et al. \cite{38} replaced the parameter $\alpha$ in (2,2)-block of the SS preconditioner by another parameter $\beta$, and employed the generalized SS (GSS) preconditioner of the form
\begin{eqnarray}
\mathcal{P}_{GSS}=\frac{1}{2}\left(
                         \begin{array}{cc}
                           \alpha I+A & B \\
                           -B^{T} & \beta I \\
                         \end{array}
                       \right),\nonumber
\end{eqnarray}
where $\alpha\geq{0}$, $\beta>0$ and $I$ is the identity matrix. It is easy to see that $\mathcal{P}_{SS}$ is a special case of $\mathcal{P}_{GSS}$ when $\alpha=\beta$. Numerical results in \cite{30,37} confirmed that the GSS preconditioner is superior to the SS preconditioner.

Very recently, based on the well-known Hermitian and skew-Hermitian splitting (HSS) of the matrix $A$: $A=H+S$, where $H=\frac{1}{2}(A+A^{T})$, $S=\frac{1}{2}(A-A^{T})$, and similar to the shift-splitting \cite{35,40}, the modified shift-splitting (MSS) preconditioner \cite{39} was proposed for nonsymmetric saddle point problem (1), the form of $\mathcal{P}_{MSS}$ is:
\begin{eqnarray}
\mathcal{P}_{MSS}=\frac{1}{2}\left(
                         \begin{array}{cc}
                           \alpha I+2H & B \\
                           -B^{T} & \alpha I \\
                         \end{array}
                       \right)\nonumber
\end{eqnarray}
with $\alpha>0$ being a constant and $I$ being the identity matrix with appropriate dimension.

In the sequel, by replacing the parameter $\alpha$ in (2,2)-block in the MSS preconditioner by another parameter $\beta$, Huang et al. \cite{59} established the generalized MSS (GMSS) preconditioner. They discussed the corresponding GMSS iteration method is convergent and semi-convergent under proper conditions, and showed that the GMSS iteration method and the GMSS preconditioner are better than the MSS iteration method and the MSS preconditioner, respectively by numerical experiments.

In order to increase the convergence rate of the GSS method for the nonsingular saddle point problems with symmetric positive definite (1,1) parts, Huang and Su \cite{51} newly developed the modified shift-splitting (MSSP) preconditioner of the form:
\begin{eqnarray}
\mathcal{P}_{MSSP}=\left(
                         \begin{array}{cc}
                           \alpha I+2A & 2B \\
                           -2B^{T} & \alpha I \\
                         \end{array}
                       \right)\nonumber
\end{eqnarray}
with $\alpha>0$ being a constant and $I$ being the identity matrix with appropriate dimension, which derived from the following modified shift-splitting of the saddle point matrix $\mathcal{A}$:
\begin{eqnarray}
\mathcal{A}=\mathcal{P}_{MSSP}-\mathcal{Q}_{MSSP}=\left(
                         \begin{array}{cc}
                           \alpha I+2A & 2B \\
                           -2B^{T} & \alpha I \\
                         \end{array}
                       \right)-\left(
                         \begin{array}{cc}
                           \alpha I+A & B \\
                           -B^{T} & \alpha I \\
                         \end{array}
                       \right).\nonumber
\end{eqnarray}
The authors in \cite{51} theoretically verified the corresponding MSSP iteration method is unconditional convergent and estimated the bounds of the eigenvalues of the iteration matrix of the MSSP iteration method. Numerical experiments illustrated that the MSSP preconditioner outperforms the SS and GSS preconditioners for the nonsingular saddle point problems with symmetric positive definite (1,1) parts.

To further accelerate the convergence rates of the GSS and the GMSS preconditioned GMRES methods for the saddle point problems with nonsymmetric positive definite (1,1) parts, a new preconditioner which is referred to as the modified generalized shift-splitting (MGSSP) preconditioner is developed for nonsymmetric saddle point problems in this paper. Theoretical analysis also shows that the corresponding splitting iteration method is convergent and semi-convergent unconditionally. Besides, we investigate the spectral properties of the corresponding preconditioned matrix and show that it has clustered eigenvalue distribution by choosing proper parameters. Numerical experiments are presented to confirm the effectiveness of the MGSSP iteration method and the MGSSP preconditioned GMRES method for solving the nonsymmetric saddle point problems.

The outline of this paper is organized as follows. In Section 2, we propose the MGSSP iteration method which induces the MGSSP preconditioner. The unconditional convergent and semi-convergent properties of the MGSSP iteration method will be proved in Sections 3 and 4, respectively. The spectral properties of the MGSSP preconditioned matrix are obtained correspondingly in Section 5. We examine the feasibility and effectiveness of the MGSSP iteration method and the MGSSP preconditioned GMRES method for solving the nonsymmetric nonsingular and singular saddle point problems by numerical experiments in Section 6. Finally, a
brief conclusion will be given to end this work in Section 7.

Throughout this paper, $\lambda_{\min}(A)$ and $\rho(A)$ represent the minimum eigenvalue and the spectral radius of the matrix $A$, respectively. $(.)^{*}$ denotes the conjugate transpose of either a vector or a matrix.

\section{The modified generalized shift-splitting (MGSSP) preconditioner and its implementation}
\label{sec:Preliminaries}

In this section, inspired by the ideas of \cite{37,38,51}, we develop a new splitting called the modified generalized shift-splitting (MGSSP) of the nonsymmetric saddle point matrix $\mathcal{A}$ by combining the generalized splitting-splitting and the modified shift-splitting of the saddle point matrix $\mathcal{A}$ as follows.
\begin{eqnarray}
\mathcal{A}=\mathcal{P}_{MGSSP}-\mathcal{Q}_{MGSSP}=\left(
                         \begin{array}{cc}
                           \alpha I+2A & 2B \\
                           -2B^{T} & \beta I \\
                         \end{array}
                       \right)-\left(
                         \begin{array}{cc}
                           \alpha I+A & B \\
                           -B^{T} & \beta I \\
                         \end{array}
                       \right),
\end{eqnarray}
where $\alpha\geq{0}$, $\beta>0$ are two constants and $I$ is the unit matrix with appropriate dimension. Then, the splitting (2) naturally leads to the following modified generalized shift-splitting iteration method for solving the nonsymmetric saddle point problem (1):

\textbf{The modified generalized shift-splitting (MGSSP) iteration method:} Let $\alpha\geq{0}$ and $\beta>0$ be two given constants. Given an initial guess $(x^{(0)^{T}},y^{(0)^{T}})^{T}$. For $k=0,1,2,\cdots$, until $(x^{(k)^{T}},y^{(k)^{T}})^{T}$ converges, compute
\begin{eqnarray}
\left(
                         \begin{array}{cc}
                           \alpha I+2A & 2B \\
                           -2B^{T} & \beta I \\
                         \end{array}
                       \right)\left(
                                \begin{array}{c}
                                  x^{(k+1)} \\
                                  y^{(k+1)} \\
                                \end{array}
                              \right)
                       =\left(
                         \begin{array}{cc}
                           \alpha I+A & B \\
                           -B^{T} & \beta I \\
                         \end{array}
                       \right)\left(
                                \begin{array}{c}
                                  x^{(k)} \\
                                  y^{(k)} \\
                                \end{array}
                              \right)+\left(
                   \begin{array}{c}
                     f \\
                     -g \\
                   \end{array}
                 \right).\nonumber
\end{eqnarray}
Hence the MGSSP iteration method can be written in the following fixed point form
\begin{eqnarray}
\left(
                                \begin{array}{c}
                                  x^{(k+1)} \\
                                  y^{(k+1)} \\
                                \end{array}
                              \right)
                       =\mathcal{T}(\alpha,\beta)\left(
                                \begin{array}{c}
                                  x^{(k)} \\
                                  y^{(k)} \\
                                \end{array}
                              \right)+\left(
                         \begin{array}{cc}
                           \alpha I+2A & 2B \\
                           -2B^{T} & \beta I \\
                         \end{array}
                       \right)^{-1}\left(
                   \begin{array}{c}
                     f \\
                     -g \\
                   \end{array}
                 \right),
\end{eqnarray}
where
\begin{eqnarray}
\mathcal{T}(\alpha,\beta)=\left(
                         \begin{array}{cc}
                           \alpha I+2A & 2B \\
                           -2B^{T} & \beta I \\
                         \end{array}
                       \right)^{-1}\left(
                         \begin{array}{cc}
                           \alpha I+A & B \\
                           -B^{T} & \beta I \\
                         \end{array}
                       \right)\nonumber
\end{eqnarray}
is the iteration matrix.

It should be noted that any matrix splitting not only can automatically lead to a splitting iteration method, but also can naturally induce a splitting preconditioner for the Krylov subspace methods. The splitting preconditioner corresponds to the MGSSP iteration (2) is given by
\begin{eqnarray}
\mathcal{P}_{MGSSP}=\left(
                         \begin{array}{cc}
                           \alpha I+2A & 2B \\
                           -2B^{T} & \beta I \\
                         \end{array}
                       \right),
\end{eqnarray}
which is called the MGSSP preconditioner for the nonsymmetric saddle point matrix $\mathcal{A}$.

At each step of the MGSSP iteration (3) or applying the MGSSP preconditioner $\mathcal{P}_{MGSSP}$ within a Krylov subspace method, we need to solve a linear system with $\mathcal{P}_{MGSSP}$ as the coefficient matrix. That is to say, we need to solve a linear system of the form
\begin{eqnarray}
\left(
                         \begin{array}{cc}
                           \alpha I+2A & 2B \\
                           -2B^{T} & \beta I \\
                         \end{array}
                       \right)z
                       =r,\nonumber
\end{eqnarray}
where $z=(z_{1}^{T},z_{2}^{T})^{T}$ and $r=(r_{1}^{T},r_{2}^{T})^{T}$ with $z_{1},r_{1}\in{\mathbb{R}}^{m}$ and $z_{2},r_{2}\in{\mathbb{R}}^{n}$. It is not difficult to check that
\begin{eqnarray}
\mathcal{P}_{MGSSP}=\left(
                                \begin{array}{cc}
                                  I & \frac{2}{\beta}B \\
                                  0 & I \\
                                \end{array}
                              \right)\left(
                                       \begin{array}{cc}
                                         \alpha I+2A+\frac{4}{\beta}BB^{T} & 0 \\
                                         0 & \beta I \\
                                       \end{array}
                                     \right)\left(
                                              \begin{array}{cc}
                                                I & 0 \\
                                                -\frac{2}{\beta}B^{T} & I \\
                                              \end{array}
                                            \right).
\end{eqnarray}
It follows from the decomposition of $\mathcal{P}_{MGSSP}$ in (5) that
\begin{eqnarray}
\left(
  \begin{array}{c}
    z_{1} \\
    z_{2} \\
  \end{array}
\right)
=\left(
                                              \begin{array}{cc}
                                                I & 0 \\
                                                \frac{2}{\beta}B^{T} & I \\
                                              \end{array}
                                            \right)\left(
                                       \begin{array}{cc}
                                         \alpha I+2A+\frac{4}{\beta}BB^{T} & 0 \\
                                         0 & \beta I \\
                                       \end{array}
                                     \right)^{-1}\left(
                                \begin{array}{cc}
                                  I & -\frac{2}{\beta}B \\
                                  0 & I \\
                                \end{array}
                              \right)\left(
                                       \begin{array}{c}
                                         r_{1} \\
                                         r_{2} \\
                                       \end{array}
                                     \right)
                              .
\end{eqnarray}
Therefore, we can derive the following algorithmic version of the MGSSP iteration method.\\
\textbf{Algorithm 2.1} For a given vector $r=(r_{1}^{T},r_{2}^{T})^{T}$, the vector $z=(z_{1}^{T},z_{2}^{T})^{T}$ can be computed by (6) according to the following steps:\\
(1)~compute $t_{1}=r_{1}-\frac{2}{\beta}Br_{2}$;\\
(2)~solve $(\alpha I+2A+\frac{4}{\beta}BB^{T})z_{1}=t_{1}$;\\
(3)~compute $z_{2}=\frac{1}{\beta}(2B^{T}z_{1}+r_{2})$.

From Algorithm 2.1, it is known that at each iteration, a linear system with the coefficient matrix $\alpha I+2A+\frac{4}{\beta}BB^{T}$ only needs to be solved. However, it may be very costly and impractical in actual implementations because of the sparsity pattern of $\alpha I+2A+\frac{4}{\beta}BB^{T}$. Fortunately, the matrix $\alpha I+2A+\frac{4}{\beta}BB^{T}$ is positive definite for all $\alpha\geq{0}$ and $\beta>0$. Therefore, we can employ the Krylov subspace method, such as the GMRES method to solve the sub-linear systems with the coefficient matrix $\alpha I+2A+\frac{4}{\beta}BB^{T}$ by a prescribed accuracy. In addition, it can be solved by some direct methods, such as the sparse LU factorization. What we want to pose here is that we always use the sparse LU factorization to solve this problem in our paper.

\section{Convergence of the MGSSP iteration method for nonsingular saddle point problems}
The main purpose of this section is to study the convergence properties of the MGSSP iteration method by
analyzing the spectral properties of the iteration matrix. Before doing this, we derive some lemmas which will be useful
in the following proofs.

\newtheorem{lem1}{Lemma}[section]
\begin{lem1}\label{lem1}\emph{\cite{5}}
Both roots of the complex quadratic equation $x^{2}-\phi x+\psi=0$ are less than one in modulus if and only if $|\phi-\bar{\phi}\psi|+|\psi|^{2}<1$, where $\bar{\phi}$ denotes the conjugate complex of $\phi$.
\end{lem1}
\newtheorem{lem4}[lem1]{Lemma}
\begin{lem4}\label{lem2}
Let $A\in{\mathbb{R}}^{m\times{m}}$ be a positive definite matrix, $B\in{\mathbb{R}}^{m\times{n}}$ be of full column rank, and $\alpha\geq{0}$ and $\beta>0$ be two given constants. If $\lambda$ is an eigenvalue of the iteration matrix $\mathcal{T}(\alpha,\beta)$, then $\lambda\neq{\pm1}$.
\end{lem4}
\noindent {\bf Proof.} Let $\lambda$ be an eigenvalue of the iteration matrix $\mathcal{T}(\alpha,\beta)$ of the MGSSP iteration method, and $(u^{*},v^{*})^{*}\in{\mathbb{C}}^{m+n}$ be the corresponding eigenvector. Then it holds that
\begin{eqnarray}
\left(
                      \begin{array}{cc}
                        \alpha I+A & B \\
                        -B^{T} & \beta I \\
                      \end{array}
                    \right)\left(
             \begin{array}{c}
               u \\
               v \\
             \end{array}
           \right)=\lambda\left(
              \begin{array}{cc}
                \alpha I+2A & 2B \\
                -2B^{T} & \beta I \\
              \end{array}
            \right)\left(
                            \begin{array}{c}
                              u \\
                              v \\
                            \end{array}
                          \right).\nonumber
\end{eqnarray}
After proper manipulations, we obtain
\begin{equation}
  \left\{
   \begin{aligned}
   &(\alpha I+A)u+Bv=\lambda(\alpha I+2A)u+2\lambda B v, \\
   &-B^{T}u+\beta v=-2\lambda B^{T}u+\lambda \beta v. \\
   \end{aligned}
   \right.
\end{equation}
Now we will give the proof by contradiction. If $\lambda=1$, then from (7), it has $Au+Bv=0$ and $B^{T}u=0$, which lead to $u=-A^{-1}Bv$ and $B^{T}A^{-1}Bv=0$. Thus we get $Bv=0$ by the positive definiteness of $A^{-1}$, and therefore $v=0$ and $u=-A^{-1}Bv=0$, a contradiction. In addition, if $\lambda=-1$, then it follows from the second equation of (7) that $v=\frac{3B^{T}u}{2\beta}$. Substituting this relation into the first equation of (7) gives $\bar{A}u=(2\alpha I+3A+\frac{9BB^{T}}{2\beta})u=0$, then $u=0$ is due to the fact that $\bar{A}$ is nonsingular, which yields that $v=\frac{3B^{T}u}{2\beta}=0$, a contradiction. \hfill$\blacksquare$

\newtheorem{lem5}[lem1]{Lemma}
\begin{lem5}\label{lem2}
Assume that the conditions in Lemma 3.2 are satisfied. Let $\lambda$ be an eigenvalue of the iteration matrix $\mathcal{T}(\alpha,\beta)$ of the MGSSP iteration method and $\mathbf{u}=(u^{*},v^{*})^{*}\in{\mathbb{C}}^{m+n}$, with $u\in{\mathbb{C}}^{m}$ and $v\in{\mathbb{C}}^{n}$, be the corresponding eigenvector. Then $u\neq0$. Moreover, if $v=0$, then $|\lambda|<1$.
\end{lem5}
\noindent {\bf Proof.} If $u=0$, then from the second equation of (7), we have
$(\lambda-1)\beta v=0$. Inasmuch as $\lambda\neq{1}$ and $\beta>0$, we derive $v=0$. This contradicts to the assumption that $\mathbf{u}=(u^{*},v^{*})^{*}$ is an eigenvector. Furthermore, if $v=0$, then it follows from the first equation of (7) that
\begin{eqnarray}
(\alpha I+A)u=\lambda(\alpha I+2A)u.
\end{eqnarray}
Since $u\neq{0}$, the definition $\frac{u^{*}}{u^{*}u}$ does make sense. Premultiplying (8) with $\frac{u^{*}}{u^{*}u}$ gives
\begin{eqnarray}
\lambda=\frac{(\alpha+a)+ib}{(\alpha+2a)+2ib},
\end{eqnarray}
where $a+ib=\frac{u^{*}Au}{u^{*}u}$.
Since $A$ is positive definite, $a>0$. It follows from (9) that
\begin{eqnarray}
|\lambda|=\sqrt{\frac{(\alpha+a)^{2}+b^{2}}{(\alpha+2a)^{2}+4b^{2}}}<1.\nonumber
\end{eqnarray}
Thus, we completes our proof of Lemma 3.3. \hfill$\blacksquare$

\newtheorem{thm1}{Theorem}[section]
\begin{thm1}\label{thm1}
Assume the conditions in Lemma 3.2 are satisfied. Let $\lambda$ be an eigenvalue of the iteration matrix $\mathcal{T}(\alpha,\beta)$ of the MGSSP iteration method and $\mathbf{u}=(u^{*},v^{*})^{*}\in{\mathbb{C}}^{m+n}$, with $u\in{\mathbb{C}}^{m}$ and $v\in{\mathbb{C}}^{n}$, be the corresponding eigenvector. Then
the MGSSP iteration method converges to the exact solution of the saddle point problem (1) for all $\alpha\geq{0}$ and $\beta>0$.
\end{thm1}
\noindent {\bf Proof.} By making use of Lemma 3.2, we have $\lambda\neq{1}$, then from the second equation of (7), it has
\begin{eqnarray}
v=\frac{(2\lambda-1)B^{T}u}{(\lambda-1)\beta},\nonumber
\end{eqnarray}
substituting it into the first equation of (7) results in
\begin{eqnarray}
\lambda^{2}(\alpha\beta I+2\beta A+4BB^{T})u-\lambda(2\alpha\beta I+3\beta A+4BB^{T})u+(\alpha\beta I+\beta A+BB^{T})u=0.
\end{eqnarray}
By making use of Lemma 3.3, it holds that $u\neq{0}$.
Denote
\begin{eqnarray}
a+ib=\frac{u^{*}Au}{u^{*}u},\ c=\frac{u^{*}BB^{T}u}{u^{*}u}\geq{0}.\nonumber
\end{eqnarray}
By multiplying $\frac{u^{*}}{u^{*}u}$ on (10) from the left, we have
\begin{eqnarray}
\lambda^{2}(\alpha\beta+2\beta a+4c+2\beta bi)-\lambda(2\alpha\beta+3\beta a+4c+3\beta bi)+(\alpha\beta+\beta a+c+\beta bi)=0.
\end{eqnarray}
Having mind that $A$ is positive definite, we get $a>0$ and $c\geq{0}$, which lead to $\alpha\beta+2\beta a+4c+2\beta bi\neq{0}$ by $\alpha\geq{0}$ and $\beta>0$. Hence, (11) can be rewritten as $\lambda^{2}-\phi \lambda+\psi=0$, where
\begin{eqnarray}
\phi=\frac{2\alpha\beta+3\beta a+4c+3\beta bi}{\alpha\beta+2\beta a+4c+2\beta bi},
\quad \psi=\frac{\alpha\beta+\beta a+c+\beta bi}{\alpha\beta+2\beta a+4c+2\beta bi}.\nonumber
\end{eqnarray}
If $c=0$, then (11) can be expressed as
\begin{eqnarray}
\lambda^{2}-\lambda\frac{2\alpha+3a+3bi}{\alpha+2a+2bi}
+\frac{\alpha+a+bi}{\alpha+2a+2bi}=0.
\end{eqnarray}
Solving the two roots of (12), we obtain
\begin{eqnarray}
\lambda=1\ \mathrm{or}\ \lambda=\frac{\alpha+a+bi}{\alpha+2a+2bi}.\nonumber
\end{eqnarray}
Lemma 3.2 implies that $\lambda\neq{1}$, then
\begin{eqnarray}
|\lambda|=\left|\frac{\alpha+a+bi}{\alpha+2a+2bi}\right|=\sqrt{\frac{(\alpha+a)^{2}+b^{2}}{(\alpha+2a)^{2}+4b^{2}}}<1.\nonumber
\end{eqnarray}
Now we turn to prove $|\lambda|<1$ under the condition $c>0$. According to Lemma 3.1, we know that $|\lambda|<1$ if and only if $|\phi-\bar{\phi}\psi|+|\psi|^{2}<1$. After some manipulations, we derive
\begin{eqnarray}
\phi-\bar{\phi}\psi=\frac{2\alpha\beta^{2}a+6\alpha\beta c+3\beta^{2}a^{2}+13\beta ac+12c^{2}+3\beta^{2}b^{2}+3\beta bci}{(\alpha\beta+2\beta a+4c)^{2}+4\beta^{2}b^{2}}\nonumber
\end{eqnarray}
and
\begin{eqnarray}
1-|\psi|^{2}=\frac{2\alpha\beta^{2}a+6\alpha\beta c+3\beta^{2}a^{2}+14\beta ac+15c^{2}+3\beta^{2}b^{2}}{(\alpha\beta+2\beta a+4c)^{2}+4\beta^{2}b^{2}}.\nonumber
\end{eqnarray}
Hence, $|\phi-\bar{\phi}\psi|+|\psi|^{2}<1$ is valid if and only if
\begin{eqnarray}
&&|2\alpha\beta^{2}a+6\alpha\beta c+3\beta^{2}a^{2}+13\beta ac+12c^{2}+3\beta^{2}b^{2}+3\beta bci|\nonumber\\
&=&\sqrt{(2\alpha\beta^{2}a+6\alpha\beta c+3\beta^{2}a^{2}+13\beta ac+12c^{2}+3\beta^{2}b^{2})^{2}+9\beta^{2} b^{2}c^{2}}\nonumber\\
&<&2\alpha\beta^{2}a+6\alpha\beta c+3\beta^{2}a^{2}+14\beta ac+15c^{2}+3\beta^{2}b^{2},\nonumber
\end{eqnarray}
which is equivalent to
\begin{eqnarray}
&&(2\alpha\beta^{2}a+6\alpha\beta c+3\beta^{2}a^{2}+13\beta ac+12c^{2}+3\beta^{2}b^{2})^{2}+9\beta^{2} b^{2}c^{2}\nonumber\\
&<&(2\alpha\beta^{2}a+6\alpha\beta c+3\beta^{2}a^{2}+14\beta ac+15c^{2}+3\beta^{2}b^{2})^{2}.
\end{eqnarray}
Since $a>0$, $c>{0}$, $b^{2}\geq{0}$, $\alpha\geq{0}$ and $\beta>0$, it holds that
\begin{eqnarray}
&&(2\alpha\beta^{2}a+6\alpha\beta c+3\beta^{2}a^{2}+14\beta ac+15c^{2}+3\beta^{2}b^{2})^{2}\nonumber\\
&=&[(2\alpha\beta^{2}a+6\alpha\beta c+3\beta^{2}a^{2}+13\beta ac+12c^{2}+3\beta^{2}b^{2}+(\beta ac+3c^{2})]^{2}\nonumber\\
&=&(2\alpha\beta^{2}a+6\alpha\beta c+3\beta^{2}a^{2}+13\beta ac+12c^{2}+3\beta^{2}b^{2})^{2}+(\beta ac+3c^{2})^{2}\nonumber\\
&&+2(2\alpha\beta^{2}a+6\alpha\beta c+3\beta^{2}a^{2}+13\beta ac+12c^{2}+3\beta^{2}b^{2})(\beta ac+3c^{2})\nonumber\\
&>&(2\alpha\beta^{2}a+6\alpha\beta c+3\beta^{2}a^{2}+13\beta ac+12c^{2}+3\beta^{2}b^{2})^{2}\nonumber\\
&&+(2\alpha\beta^{2}a+6\alpha\beta c+3\beta^{2}a^{2}+13\beta ac+12c^{2}+3\beta^{2}b^{2})(\beta ac+3c^{2})\nonumber\\
&>&(2\alpha\beta^{2}a+6\alpha\beta c+3\beta^{2}a^{2}+13\beta ac+12c^{2}+3\beta^{2}b^{2})^{2}+3\beta^{2}b^{2}(\beta ac+3c^{2})\nonumber\\
&\geq&(2\alpha\beta^{2}a+6\alpha\beta c+3\beta^{2}a^{2}+13\beta ac+12c^{2}+3\beta^{2}b^{2})^{2}+9\beta^{2}b^{2}c^{2},\nonumber
\end{eqnarray}
which implies that (13) holds true, i.e., $|\phi-\bar{\phi}\psi|+|\psi|^{2}<1$ and therefore $|\lambda|<1$. Hence, the MGSSP iteration method is convergent for any $\alpha\geq{0}$ and $\beta>0$. This proof is completed. \hfill$\blacksquare$

\section{Semi-convergence of the MGSSP iteration method for singular saddle point problems}
When the saddle point matrix $\mathcal{A}$ is nonsingular, the MGSSP iteration scheme (3) converges to the exact solution
of (1) for any initial vector if and only if $\rho(\mathcal{T}(\alpha,\beta))<1$, whereas for the singular matrix $\mathcal{A}$, we have $\rho(\mathcal{T}(\alpha,\beta))\geq 1$. In this section, we assume that the sub-matrix $B$ in (1) is rank deficient and discuss the semi-convergence of the MGSSP iteration method for solving the singular saddle point problems.

To analyze the semi-convergent properties of the MGSSP iteration method, we present the following lemma which describes the semi-convergence property about the iteration scheme (3) when $\mathcal{A}$ is singular.
\newtheorem{lem7}[lem1]{Lemma}
\begin{lem7}\emph{\cite{32}}
The iteration scheme (3) is semi-convergent if and only if the following two conditions are satisfied:\\
(i)~$index(I-T)=1$, or equivalently, $rank((I-T)^{2})=rank(I-T)$, where $T=I-GM$ is the iteration matrix;\\
(ii)~the pseudo-spectral radius of $T$ is less than $1$, i.e.,
\begin{eqnarray}
\gamma(T)=\max\{|\lambda|:\lambda\in\sigma(T),\lambda\neq{1}\}<1,\nonumber
\end{eqnarray}
where $\sigma(T)$ is the spectral set of the matrix $T$. Here, we denote the null space, the index and the rank of $A$ by $null(A)$, $index(A)$ and $rank(A)$, respectively.
\end{lem7}

Lemma 4.1 describes the semi-convergence property about the iteration scheme (3) when $\mathcal{A}$ is singular. Therefore, to
get the semi-convergence property of the MGSSP iteration method, only the two conditions in Lemma 4.1 need to verify. We
consider these two conditions in Lemmas 4.2 and 4.3, respectively.
\newtheorem{lem9}[lem1]{Lemma}
\begin{lem9}\label{lem9}
Let $A$ be nonsymmetric positive definite, $B$ be rank deficient and $\alpha\geq{0},\beta>0$ be given constants. Then, the iteration matrix $\mathcal{T}(\alpha,\beta)$ of the MGSSP iteration method satisfies $index(I-\mathcal{T}(\alpha,\beta))=1$, or equivalent
\begin{eqnarray}
rank(I-\mathcal{T}(\alpha,\beta))=rank((I-\mathcal{T}(\alpha,\beta))^{2}),
\end{eqnarray}
where $\mathcal{T}(\alpha,\beta)$ is the iteration matrix of the MGSSP iteration method defined as in (3).
\end{lem9}
\noindent {\bf Proof.} Inasmuch as $\mathcal{T}(\alpha,\beta)=\mathcal{P}_{MGSSP}^{-1}\mathcal{Q}_{MGSSP}=I-\mathcal{P}_{MGSSP}^{-1}\mathcal{A}$, Equation (14) holds if
\begin{eqnarray}
null(\mathcal{P}_{MGSSP}^{-1}\mathcal{A})=null((\mathcal{P}_{MGSSP}^{-1}\mathcal{A})^{2}).\nonumber
\end{eqnarray}
It is easy to see that $null(\mathcal{P}_{MGSSP}^{-1}\mathcal{A})\subseteq null((\mathcal{P}_{MGSSP}^{-1}\mathcal{A})^{2})$. Thus we only need to prove
\begin{eqnarray}
null(\mathcal{P}_{MGSSP}^{-1}\mathcal{A})\supseteq null((\mathcal{P}_{MGSSP}^{-1}\mathcal{A})^{2}).\nonumber
\end{eqnarray}
Let $x=(x_{1}^{*},x_{2}^{*})^{*}\in{\mathbb{C}^{m+n}}\in null((\mathcal{P}_{MGSSP}^{-1}\mathcal{A})^{2})$, then it has $(\mathcal{P}_{MGSSP}^{-1}\mathcal{A})^{2}x=0$. Denote by $y=\mathcal{P}_{MGSSP}^{-1}\mathcal{A}x$. After suitable manipulations, we have
\setlength\arraycolsep{3pt}
\begin{eqnarray}
y=\left(
    \begin{array}{c}
      y_{1} \\
      y_{2} \\
    \end{array}
  \right)&=&\left(
                         \begin{array}{cc}
                           \alpha I+2A & 2B \\
                           -2B^{T} & \beta I \\
                         \end{array}
                       \right)^{-1}\left(
    \begin{array}{cc}
      A & B \\
      -B^{T} & 0 \\
    \end{array}
  \right)\left(
           \begin{array}{c}
             x_{1} \\
             x_{2} \\
           \end{array}
         \right)\nonumber\\
&=&\left(
                                              \begin{array}{cc}
                                                I & 0 \\
                                                \frac{2}{\beta}B^{T} & I \\
                                              \end{array}
                                            \right)\left(
                                       \begin{array}{cc}
                                         \alpha I+2A+\frac{4}{\beta}BB^{T} & 0 \\
                                         0 & \beta I \\
                                       \end{array}
                                     \right)^{-1}\left(
                                \begin{array}{cc}
                                  I & -\frac{2}{\beta}B \\
                                  0 & I \\
                                \end{array}
                              \right)\left(
    \begin{array}{cc}
      A & B \\
      -B^{T} & 0 \\
    \end{array}
  \right)\left(
           \begin{array}{c}
             x_{1} \\
             x_{2} \\
           \end{array}
         \right)\nonumber\\
&=&\left(
   \begin{array}{c}
     \left(\alpha I+2A+\frac{4}{\beta}BB^{T}\right)^{-1}\left(Ax_{1}+Bx_{2}+\frac{2}{\beta}BB^{T}x_{1}\right) \\
     \frac{2}{\beta}B^{T}\left(\alpha I+2A+\frac{4}{\beta}BB^{T}\right)^{-1}\left(Ax_{1}+Bx_{2}+\frac{2}{\beta}BB^{T}x_{1}\right)-\frac{1}{\beta}B^{T}x_{1} \\
   \end{array}
 \right),\nonumber
\end{eqnarray}
i.e.,
\begin{equation}
  \left\{
   \begin{aligned}
   &y_{1}=\left(\alpha I+2A+\frac{4}{\beta}BB^{T}\right)^{-1}\left(Ax_{1}+Bx_{2}+\frac{2}{\beta}BB^{T}x_{1}\right)
   ,\\
   &y_{2}=\frac{2}{\beta}B^{T}\left(\alpha I+2A+\frac{4}{\beta}BB^{T}\right)^{-1}\left(Ax_{1}+Bx_{2}+\frac{2}{\beta}BB^{T}x_{1}\right)-\frac{1}{\beta}B^{T}x_{1}. \\
   \end{aligned}
   \right.
\end{equation}
Since $\mathcal{P}_{MGSSP}^{-1}\mathcal{A}y=(\mathcal{P}_{MGSSP}^{-1}\mathcal{A})^{2}x=0$, it holds that $\mathcal{A}y=0$, i.e.,
\begin{eqnarray}
 Ay_{1}+By_{2}=0,\ -B^{T}y_{1}=0.
\end{eqnarray}
Since $A$ is positive definite, from the first equation of (16) we can easily get $y_{1}=-A^{-1}By_{2}$. Then substituting this relationship into the second equation of (16), we obtain $B^{T}A^{-1}By_{2}=0$, which leads to $By_{2}=0$. Taking $By_{2}=0$ into $y_{1}=-A^{-1}By_{2}$, we obtain $y_{1}=0$. Hence, the first equation of (15) becomes
\begin{eqnarray}
y_{1}=\left(\alpha I+2A+\frac{4}{\beta}BB^{T}\right)^{-1}\left(Ax_{1}+Bx_{2}+\frac{2}{\beta}BB^{T}x_{1}\right)=0.\nonumber
\end{eqnarray}
Substituting $y_{1}=0$ into $y_{2}$ yields $y_{2}=-\frac{1}{\beta}B^{T}x_{1}$. Since $By_{2}=0$, $-\frac{1}{\beta}BB^{T}x_{1}=0$, it has $x_{1}^{*}BB^{T}x_{1}=0$. This results in $B^{T}x_{1}=0$, then we get $y_{2}=-\frac{1}{\beta}B^{T}x_{1}=0$. Thus, $y=\mathcal{P}_{MGSSP}^{-1}\mathcal{A}x=0$, i.e.,
\begin{eqnarray}
null(\mathcal{P}_{MGSSP}^{-1}\mathcal{A})\supseteq null((\mathcal{P}_{MGSSP}^{-1}\mathcal{A})^{2}).
\end{eqnarray}
The conclusion follows by (17). \hfill$\blacksquare$

In the sequel, we show that the iteration scheme (3) satisfies the condition (ii) in Lemma 4.1. Let $B=U(B_{r},0)V^{*}$ be the singular decomposition of matrix $B$, where
\begin{eqnarray}
B_{r}=\left(
        \begin{array}{c}
          \Sigma_{r} \\
          0 \\
        \end{array}
      \right)\in{\mathbb{C}^{m\times{r}}},\quad \Sigma_{r}=\mathrm{diag}(\sigma_{1},\sigma_{2},\cdots,\sigma_{r})\in{\mathbb{C}^{r\times{r}}}\nonumber
\end{eqnarray}
with $U\in{\mathbb{C}^{m\times{m}}}$, $V\in{\mathbb{C}^{n\times{n}}}$ being two unitary matrices and $\sigma_{i}$ $(i=1,2,\cdots,r)$ being a singular value of $B$.

We introduce a block diagonal matrix
\setlength\arraycolsep{4pt}
\begin{eqnarray}
P=\left(
    \begin{array}{cc}
      U & 0 \\
      0 & V \\
    \end{array}
  \right)\nonumber
\end{eqnarray}
which is a $(m+n)\times{(m+n)}$ unitary matrix, and the iteration matrix $\mathcal{T}(\alpha,\beta)$ is unitary similar to the matrix $\hat{\mathcal{T}}(\alpha,\beta)=P^{*}\mathcal{T}(\alpha,\beta)P$. Hence, the matrix $\mathcal{T}(\alpha,\beta)$ has the same spectrum with the matrix $\hat{\mathcal{T}}(\alpha,\beta)$. Thus we only need to analyze the pseudo-spectral radius of the matrix $\hat{\mathcal{T}}(\alpha,\beta)$ now.

Denoting $\hat{A}=U^{*}AU$, then it holds that
\begin{eqnarray}
\hat{\mathcal{T}}(\alpha,\beta)&=&P^{*}\left(
                         \begin{array}{cc}
                           \alpha I+2A & 2B \\
                           -2B^{T} & \beta I \\
                         \end{array}
                       \right)^{-1}\left(
                         \begin{array}{cc}
                           \alpha I+A & B \\
                           -B^{T} & \beta I \\
                         \end{array}
                       \right)P\nonumber\\
&=&\left(\begin{array}{cc}
                           \alpha I+2U^{*}AU & 2U^{*}BV \\
                           -2V^{*}B^{T}U & \beta I \\
                         \end{array}
                       \right)^{-1}\left(
                         \begin{array}{cc}
                           \alpha I+U^{*}AU & U^{*}BV \\
                           -V^{*}B^{T}U & \beta I \\
                         \end{array}
                       \right)\nonumber\\
&=&\left(
     \begin{array}{ccc}
       \alpha I+2\hat{A} & 2B_{r} & 0 \\
       -2B_{r}^{T} & \beta I & 0 \\
       0 & 0 & \beta I \\
     \end{array}
   \right)^{-1}\left(
            \begin{array}{ccc}
              \alpha I+\hat{A} & B_{r} & 0 \\
              -B_{r}^{T} & \beta I & 0 \\
              0 & 0 & \beta I \\
            \end{array}
          \right)\nonumber\\
&=&\left(
     \begin{array}{cc}
       \left(
          \begin{array}{cc}
            \alpha I+2\hat{A} & 2B_{r} \\
            -2B_{r}^{T} & \beta I \\
          \end{array}
        \right)^{-1}\left(
                      \begin{array}{cc}
                        \alpha I+\hat{A} & B_{r} \\
                        -B_{r} & \beta I \\
                      \end{array}
                    \right)
        & 0 \\
       0 & I_{n-r} \\
     \end{array}
   \right)\nonumber\\
&=&\left(
     \begin{array}{cc}
       \tilde{\mathcal{T}}(\alpha,\beta) & 0 \\
       0 & I_{n-r} \\
     \end{array}
   \right).
\end{eqnarray}
Then, from Equation (18), $\gamma(\hat{\mathcal{T}}(\alpha,\beta))<1$ holds if and only if $\rho(\tilde{\mathcal{T}}(\alpha,\beta))<1$.

Note that $\tilde{\mathcal{T}}(\alpha,\beta)$ can be viewed as the iteration matrix of the MGSSP iteration method applied to the nonsymmetric nonsingular saddle point problem
\begin{eqnarray}
\left(
    \begin{array}{cc}
      \hat{A} & B_{r} \\
      -B_{r}^{T} & 0 \\
    \end{array}
  \right)\left(
           \begin{array}{c}
             \hat{x} \\
             \hat{y} \\
           \end{array}
         \right)=\left(
                   \begin{array}{c}
                     \hat{f} \\
                     -\hat{g} \\
                   \end{array}
                 \right),\nonumber
\end{eqnarray}
where $\hat{A}=U^{*}AU$ and $\hat{y},\hat{g}\in{\mathbb{R}^{r}}$.

$\rho(\tilde{\mathcal{T}}(\alpha,\beta))<1$ implies $\gamma(\mathcal{T}(\alpha,\beta))=\gamma(\hat{\mathcal{T}}(\alpha,\beta))<1$. By making use of the proof of Theorem 3.1, we derive the following result.
\newtheorem{lem10}[lem1]{Lemma}
\begin{lem10}\label{lem10}
Let $A$ be nonsymmetric positive definite, $B$ be rank deficient and $\alpha\geq{0},\beta>0$ be two given constants. Then, the pseudo-spectral radius of the matrix $\gamma(\alpha,\beta)$ is less than 1, i.e., $\mathcal{V}(\mathcal{T}(\alpha,\beta))<1$ for all $\alpha\geq{0}$ and $\beta>0$.
\end{lem10}
It follows from Lemmas 4.2 and 4.3 that two conditions in Lemma 4.1 are satisfied. Thus, the semi-convergence of the MGSSP iteration method for solving nonsymmetric singular saddle point problems can be obtained in the following theorem.
\newtheorem{thm2}[thm1]{Theorem}
\begin{thm2}\label{thm1}
Let $A$ be nonsymmetric positive definite, $B$ be rank deficient and $\alpha\geq{0},\beta>0$ be two given constants. Then the MGSSP iteration method is semi-convergent for solving the nonsymmetric singular saddle point problem (1) for all $\alpha\geq{0}$ and $\beta>0$.
\end{thm2}

\section{Spectral analysis of the MGSSP preconditioned matrix}
The MGSSP iteration method is a stationary iteration method. Although the unconditional convergence and semi-convergence properties of the MGSSP iteration method are studied in Theorem 3.1 and Theorem 4.1, respectively, the convergence (semi-convergence) rates of the MGSSP iteration method may be slow even with the optimal parameters. To accelerate the convergence (semi-convergence) rates of the MGSSP iteration method, we consider applying the preconditioning techniques. In general, the eigenvalue and eigenvector distributions of the preconditioned matrix relate closely to the convergence
rates of Krylov subspace methods. Therefore, it is of significance to investigate the spectral properties of the preconditioned
matrix $\mathcal{P}_{MGSSP}^{-1}\mathcal{A}$. In this section, some spectral properties of the preconditioned matrix $\mathcal{P}_{MGSSP}^{-1}\mathcal{A}$ are studied.

\newtheorem{thm4}[thm1]{Theorem}
\begin{thm4}\label{thm4}
Let the MGSSP preconditioner be defined as in (4) and $(\lambda,(u^{*},v^{*})^{*})$ be an eigenpair of the preconditioned matrix $\mathcal{P}_{MGSSP}^{-1}\mathcal{A}$. Then if $B$ is of full column rank and $B^{T}u=0$, then
\begin{eqnarray}
\frac{\lambda_{\min}(H)(\alpha+2\lambda_{\min}(H))}{(\alpha+2\rho(H))^{2}+4\rho(S)^{2}}\leq Re(\lambda)
\leq \frac{\rho(H)(\alpha+2\rho(H))+2\rho(S)^{2}}{(\alpha+2\lambda_{\min}(H))^{2}},\quad |Im(\lambda)|\leq\frac{\alpha\rho(S)}{(\alpha+2\lambda_{\min}(H))^{2}},
\end{eqnarray}
where $Re(\lambda)$ and $Im(\lambda)$ denote the real part and the imaginary part of $\lambda$, respectively. If $B$ is rank deficient and $u=0$, then $\lambda=0$. Besides, if $B$ is rank deficient and $B^{T}u=0$, then $\lambda=0$ or $\lambda$ satisfies the Inequalities (19). If $B^{T}u\neq0$, then the eigenvalues of the preconditioned matrix $\mathcal{P}_{MGSSP}^{-1}\mathcal{A}$ satisfy
\begin{eqnarray}
\lambda_{+}=\frac{1}{2}+\frac{(z_{1}-\alpha\beta-\beta a_{1})+i(z_{2}-\beta b_{1})}{2(\alpha\beta+2\beta a_{1}+4c_{1}+2i\beta b_{1})},\ \lambda_{-}=\frac{1}{2}-\frac{(z_{1}+\alpha\beta+\beta a_{1})+i(z_{2}+\beta b_{1})}{2(\alpha\beta+2\beta a_{1}+4c_{1}+2i\beta b_{1})},
\end{eqnarray}
where
\begin{eqnarray}
\frac{u^{*}Au}{u^{*}u}=a_{1}+ib_{1},\ \frac{u^{*}BB^{T}u}{u^{*}u}=c_{1}
\end{eqnarray}
and $z_{1},z_{2}$ are real numbers and $z_{1}+iz_{2}$ is one of the square roots of $a_{2}+b_{2}i$, with
\begin{eqnarray}
a_{2}=\beta^{2}(a_{1}^{2}-b_{1}^{2})-4\alpha\beta c_{1},\ b_{2}=2\beta^{2}a_{1}b_{1}\nonumber
\end{eqnarray}
and
\begin{eqnarray}
&&z_{1}=\sqrt{\frac{\sqrt{[\beta^{2}(a_{1}^{2}-b_{1}^{2})-4\alpha\beta c_{1}]^{2}+4\beta^{4}a_{1}^{2}b_{1}^{2}
}+\beta^{2}(a_{1}^{2}-b_{1}^{2})-4\alpha\beta c_{1}}{2}},\nonumber\\
&&z_{2}=sign(b_{1})\sqrt{\frac{\sqrt{[\beta^{2}(a_{1}^{2}-b_{1}^{2})-4\alpha\beta c_{1}]^{2}+4\beta^{4}a_{1}^{2}b_{1}^{2}
}-\beta^{2}(a_{1}^{2}-b_{1}^{2})+4\alpha\beta c_{1}}{2}},
\end{eqnarray}
and the second root of $a_{2}+b_{2}i$ is $-(z_{1}+iz_{2})$. The eigenvalues $\lambda_{\pm}$ satisfy the following inequality:
\begin{eqnarray}
\left|\lambda_{\pm}-\frac{1}{2}\right|^{2}\leq f(a_{1},b_{1},c_{1})\leq
\frac{(\alpha\beta+2\beta \rho(H)^{2})^{2}+(\beta\rho(S)+\sqrt{\beta^{2}\rho(S)^{2}+4\alpha\beta \rho(BB^{T})})^{2}}{4(\alpha\beta+2\beta\lambda_{\min}(H)+4\lambda_{\min}(BB^{T}))^{2}}.
\end{eqnarray}
When $\beta\rightarrow{0_{+}}$, it holds that
\begin{equation}
 \left\{
   \begin{aligned}
   &\lambda_{+}=\frac{1}{2}+\frac{(z_{1}-\alpha\beta-\beta a_{1})+i(z_{2}-\beta b_{1})}{2(\alpha\beta+2\beta a_{1}+4c_{1}+2i\beta b_{1})}\rightarrow \frac{1}{2},\\
   &\lambda_{-}=\frac{1}{2}-\frac{(z_{1}+\alpha\beta+\beta a_{1})+i(z_{2}+\beta b_{1})}{2(\alpha\beta+2\beta a_{1}+4c_{1}+2i\beta b_{1})}\rightarrow \frac{1}{2},\\
   \end{aligned}
   \right.\nonumber
\end{equation}
i.e., for $\alpha>0$, the eigenvalues of the preconditioned matrix $\mathcal{P}_{MGSSP}^{-1}\mathcal{A}$ tend to scatter near the point $(\frac{1}{2},0)$ as $\beta\rightarrow{0_{+}}$; and when $\alpha\rightarrow{0_{+}}$, it has
\begin{equation}
 \left\{
   \begin{aligned}
   &\lambda_{+}=\frac{1}{2}+\frac{(z_{1}-\alpha\beta-\beta a_{1})+i(z_{2}-\beta b_{1})}{2(\alpha\beta+2\beta a_{1}+4c_{1}+2i\beta b_{1})}\rightarrow \frac{1}{2}+\frac{(z_{1}-\beta a_{1})+i(z_{2}-\beta b_{1})}{2(2\beta a_{1}+4c_{1}+2i\beta b_{1})}=\frac{1}{2},\\
   &\lambda_{-}=\frac{1}{2}-\frac{(z_{1}+\alpha\beta+\beta a_{1})+i(z_{2}+\beta b_{1})}{2(\alpha\beta+2\beta a_{1}+4c_{1}+2i\beta b_{1})}\rightarrow \frac{1}{2}-\frac{(z_{1}+\beta a_{1})+i(z_{2}+\beta b_{1})}{2(2\beta a_{1}+4c_{1}+2i\beta b_{1})}
   =\frac{1}{2}-\frac{\beta a_{1}+i\beta b_{1}}{2\beta a_{1}+4c_{1}+2i\beta b_{1}}.\\
   \end{aligned}
   \right.\nonumber
\end{equation}
That is, for $\beta>0$, the eigenvalues of the preconditioned matrix $\mathcal{P}_{MGSSP}^{-1}\mathcal{A}$ tend to scatter near the point $(\frac{1}{2},0)$ and the point $(\frac{\beta a_{1}c_{1}+2c_{1}^{2}}{(\beta a_{1}^{2}+2c_{1})^{2}+\beta^{2}b_{1}^{2}},-\frac{\beta b_{1}c_{1}}{(\beta a_{1}^{2}+2c_{1})^{2}+\beta^{2}b_{1}^{2}})$ as $\alpha\rightarrow{0_{+}}$.

In addition, the eigenvalues of $\mathcal{P}_{MGSSP}^{-1}\mathcal{A}$ tend to scatter near the points \\ $(\frac{\alpha_{0}\beta_{0}^{2}a_{1}+2\beta_{0}^{2}(a_{1}^{2}+b_{1}^{2})+12\beta_{0} a_{1}c_{1}+(\alpha_{0}\beta_{0}+4c_{1})(4c_{1}+z_{1})+2\beta_{0}(a_{1}z_{1}+|b_{1}z_{2}|)}{2[(\alpha_{0}\beta_{0}+2\beta_{0}a_{1}+4c_{1})^{2}+4\beta_{0}^{2}b_{1}^{2}]}$,  $\frac{(\alpha_{0}\beta_{0}+2\beta_{0} a_{1}+4c_{1})(z_{2}+\beta_{0} b_{1})-2\beta_{0} b_{1}(\beta_{0} a_{1}+4c_{1}+z_{1})}{2[(\alpha_{0}\beta_{0}+2\beta_{0}a_{1}+4c_{1})^{2}+4\beta_{0}^{2}b_{1}^{2}]})$ and $(\frac{\alpha_{0}\beta_{0}^{2}a_{1}+2\beta_{0}^{2}(a_{1}^{2}+b_{1}^{2})+12\beta_{0} a_{1}c_{1}+(\alpha_{0}\beta_{0}+4c_{1})(4c_{1}-z_{1})-2\beta_{0}(a_{1}z_{1}+|b_{1}z_{2}|)}{2[(\alpha_{0}\beta_{0}+2\beta_{0}a_{1}+4c_{1})^{2}+4\beta_{0}^{2}b_{1}^{2}]},
\frac{(\alpha_{0}\beta_{0}+2\beta_{0} a_{1}+4c_{1})(\beta_{0} b_{1}-z_{2})-2\beta_{0} b_{1}(\beta_{0} a_{1}+4c_{1}-z_{1})}{2[(\alpha_{0}\beta_{0}+2\beta_{0}a_{1}+4c_{1})^{2}+4\beta_{0}^{2}b_{1}^{2}]})$ as $\alpha\rightarrow \alpha_{0}$ and $\beta\rightarrow\beta_{0}$ $(0\leq\alpha_{0}<+\infty,0<\beta_{0}<+\infty)$.
\end{thm4}
\noindent {\bf Proof.} Let $(\lambda,(u^{*},v^{*})^{*})$ be an eigenpair of the preconditioned matrix $\mathcal{P}_{MGSSP}^{-1}\mathcal{A}$, we consider the eigenvalue problem $\mathcal{P}_{MGSSP}^{-1}\mathcal{A}\eta=\lambda\eta$, where $\eta=(u^{*},v^{*})^{*}$, then it holds that
\begin{eqnarray}
\left(
                         \begin{array}{cc}
                           A & B \\
                           -B^{T} & 0 \\
                         \end{array}
                       \right)\left(
                                \begin{array}{c}
                                  u \\
                                  v \\
                                \end{array}
                              \right)
                       =\lambda\left(
                         \begin{array}{cc}
                           \alpha I+2A & 2B \\
                           -2B^{T} & \beta I \\
                         \end{array}
                       \right)\left(
                                \begin{array}{c}
                                  u \\
                                  v \\
                                \end{array}
                              \right).\nonumber
\end{eqnarray}
By simple manipulations, we get
\begin{equation}
 \left\{
   \begin{aligned}
   &Au+Bv=\lambda(\alpha I+2A)u+2\lambda Bv,\\
   &-B^{T}u=-2\lambda B^{T}u+\lambda\beta v,\\
   \end{aligned}
   \right.\nonumber
\end{equation}
i.e.,
\begin{equation}
 \left\{
   \begin{aligned}
   &Au=\lambda(\alpha I+2A)u+(2\lambda-1) Bv,\\
   &(2\lambda-1)B^{T}u=\lambda\beta v.\\
   \end{aligned}
   \right.
\end{equation}
If $B$ has full column rank and $u=0$, then it follows from the second equation of (24) that $\lambda v=0$ and therefore $v=0$, which contradicts to the assumption that $(u^{*},v^{*})^{*}$ is an eigenvector. Hence $u\neq{0}$. If $B$ is of full column rank and $B^{T}u=0$, then from the second equation of (24), we have $v=0$ and
\begin{eqnarray}
Au=\lambda(\alpha I+2A)u.
\end{eqnarray}
Owing to $u\neq0$, it holds that the definition $\frac{u^{*}}{u^{*}u}$ does make sense. Premultiplying Equation (25) with $\frac{u^{*}}{u^{*}u}$ and utilizing the symbols defined as in (21) give
\begin{eqnarray}
\lambda=\frac{a_{1}+ib_{1}}{\alpha+2a_{1}+2ib_{1}}=\frac{a_{1}(\alpha+2a_{1})+2b_{1}^{2}+i\alpha b_{1}}{(\alpha+2a_{1})^{2}+4b_{1}^{2}}.
\end{eqnarray}
It is easy to verify that $\lambda\rightarrow\frac{1}{2}$ as $\alpha\rightarrow 0_{+}$. Besides, (26) implies that
\begin{eqnarray}
Re(\lambda)=\frac{a_{1}(\alpha+2a_{1})+2b_{1}^{2}}{(\alpha+2a_{1})^{2}+4b_{1}^{2}},\quad Im(\lambda)=\frac{\alpha b_{1}}{(\alpha+2a_{1})^{2}+4b_{1}^{2}}.\nonumber
\end{eqnarray}
Since
\begin{eqnarray}
&&\lambda_{\min}(H)\leq a_{1}=\frac{1}{2}\left(\frac{u^{*}Au}{u^{*}u}+\frac{u^{*}A^{T}u}{u^{*}u}\right)=\frac{u^{*}Hu}{u^{*}u}=a_{1}\leq \rho(H),\nonumber\\
&&0\leq|b_{1}|=\frac{1}{2}\left|\frac{1}{i}\left(\frac{u^{*}Au}{u^{*}u}-\frac{u^{*}A^{T}u}{u^{*}u}\right)\right|=\left|\frac{u^{*}iSu}{u^{*}u}\right|\leq\rho(S),\nonumber
\end{eqnarray}
it is not difficult to derive (19).

If $B$ is rank deficient and $u=0$, then from the second equation of (24), we derive $\lambda=0$. Additionally, if $B$ is rank deficient and $B^{T}u=0$, then it holds that $\lambda=0$ or $v=0$, $\lambda\neq{0}$ by virtue of the second equation of (24). Similar to the derivation of (19), we also deduce (19) as $B$ is rank deficient, $v=0$ and $\lambda\neq{0}$.

In the sequel, we assume that $B^{T}u\neq{0}$. Then $\lambda\neq{0}$ and $u\neq{0}$. Otherwise, it follows from the second equation of (24) that $B^{T}u={0}$, a contradiction. From the second equation of (24) we can easily get $v=\frac{(2\lambda-1)B^{T}u}{\lambda\beta}$. Then substituting
this relationship into the first equation of (24) gives
\begin{eqnarray}
\lambda^{2}(\alpha\beta I+2\beta A+4BB^{T})u-\lambda(4BB^{T}+\beta A)u+BB^{T}u=0.
\end{eqnarray}
Multiplying $\frac{u^{*}}{u^{*}u}$ on Equation (27) from the left and utilizing the symbols defined as in (21) give
\begin{eqnarray}
\lambda^{2}(\alpha\beta +2\beta a_{1}+2i\beta b_{1}+4c_{1})-\lambda(4c_{1}+\beta a_{1}+i\beta b_{1})+c_{1}=0,\nonumber
\end{eqnarray}
which can be equivalently transformed into the following equation
\begin{eqnarray}
\lambda^{2}-\lambda\frac{4c_{1}+\beta a_{1}+i\beta b_{1}}{\alpha\beta +2\beta a_{1}+2i\beta b_{1}+4c_{1}}+\frac{c_{1}}{\alpha\beta +2\beta a_{1}+2i\beta b_{1}+4c_{1}}=0.
\end{eqnarray}
By solving Equation (28), we obtain its two roots as follows:
\begin{eqnarray}
\lambda_{+}=\frac{1}{2}+\frac{(z_{1}-\alpha\beta-\beta a_{1})+i(z_{2}-\beta b_{1})}{2(\alpha\beta+2\beta a_{1}+4c_{1}+2i\beta b_{1})},\ \lambda_{-}=\frac{1}{2}-\frac{(z_{1}+\alpha\beta+\beta a_{1})+i(z_{2}+\beta b_{1})}{2(\alpha\beta+2\beta a_{1}+4c_{1}+2i\beta b_{1})},
\end{eqnarray}
where $z_{1}$ and $z_{2}$ are given by (22).
Applying (22) leads to
\begin{eqnarray}
z_{1}&=&\sqrt{\frac{\sqrt{[\beta^{2}(a_{1}^{2}-b_{1}^{2})-4\alpha\beta c_{1}]^{2}+4\beta^{4}a_{1}^{2}b_{1}^{2}
}+\beta^{2}(a_{1}^{2}-b_{1}^{2})-4\alpha\beta c_{1}}{2}},\nonumber\\
&=&\sqrt{\frac{\sqrt{\beta^{4}(a_{1}^{2}+b_{1}^{2})^{2}-8\alpha c_{1}\beta^{3}(a_{1}^{2}-b_{1}^{2})+16\alpha^{2}\beta^{2}c_{1}^{2}
}+\beta^{2}(a_{1}^{2}-b_{1}^{2})-4\alpha\beta c_{1}}{2}},\nonumber\\
&\leq&\sqrt{\frac{\sqrt{[\beta^{2}(a_{1}^{2}+b_{1}^{2})+4\alpha\beta c_{1}]^{2}
}+\beta^{2}(a_{1}^{2}-b_{1}^{2})-4\alpha\beta c_{1}}{2}}=\beta a_{1},\\
|z_{2}|&=&\sqrt{\frac{\sqrt{[\beta^{2}(a_{1}^{2}-b_{1}^{2})-4\alpha\beta c_{1}]^{2}+4\beta^{4}a_{1}^{2}b_{1}^{2}
}-\beta^{2}(a_{1}^{2}-b_{1}^{2})+4\alpha\beta c_{1}}{2}},\nonumber\\
&\leq&\sqrt{\frac{\sqrt{[\beta^{2}(a_{1}^{2}+b_{1}^{2})+4\alpha\beta c_{1}]^{2}
}-\beta^{2}(a_{1}^{2}-b_{1}^{2})+4\alpha\beta c_{1}}{2}}=\sqrt{\beta^{2}b_{1}^{2}+4\alpha\beta c_{1}},
\end{eqnarray}
which yields that
\begin{eqnarray}
\left|\lambda_{\pm}-\frac{1}{2}\right|^{2}&=&\frac{(\alpha\beta+\beta a_{1}\pm z_{1})^{2}+(\beta b_{1}\pm z_{2})^{2}}{4[(\alpha\beta+2\beta a_{1}+4c_{1})^{2}+4\beta^{2}b_{1}^{2}]}\nonumber\\
&\leq&\frac{(\alpha\beta+2\beta a_{1})^{2}+(\beta|b_{1}|+\sqrt{\beta^{2}b_{1}^{2}+4\alpha\beta c_{1}})^{2}}{4[(\alpha\beta+2\beta a_{1}+4c_{1})^{2}+4\beta^{2}b_{1}^{2}]}
:=f(a_{1},b_{1},c_{1}).
\end{eqnarray}
It is evident that an upper bound of $\left|\lambda_{\pm}-\frac{1}{2}\right|^{2}$ is $f(a_{1},b_{1},c_{1})$, with $a_{1},b_{1},c_{1}$ being bounded as follows:
\begin{eqnarray}
\lambda_{\min}(H)\leq a_{1}\leq \rho(H),\ 0\leq|b_{1}|\leq\rho(S),\ 0\leq b_{1}^{2}\leq \rho(S)^{2},\ \lambda_{\min}(BB^{T})\leq c_{1}\leq \rho(BB^{T}),\nonumber
\end{eqnarray}
which leads to
\begin{eqnarray}
\left|\lambda_{\pm}-\frac{1}{2}\right|^{2}\leq f(a_{1},b_{1},c_{1})\leq
\frac{(\alpha\beta+2\beta \rho(H)^{2})^{2}+(\beta\rho(S)+\sqrt{\beta^{2}\rho(S)^{2}+4\alpha\beta \rho(BB^{T})})^{2}}{4(\alpha\beta+2\beta\lambda_{\min}(H)+4\lambda_{\min}(BB^{T}))^{2}}.\nonumber
\end{eqnarray}
Furthermore, it is not difficult to verify that $z_{1},z_{2}\rightarrow0$ as $\beta\rightarrow{0_{+}}$, and therefore for $\alpha>0$, $\lambda_{+},\lambda_{-}\rightarrow\frac{1}{2}$ as $\beta\rightarrow{0_{+}}$. Moreover, if $\alpha\rightarrow{0_{+}}$, then it follows from (22) that $z_{1}\rightarrow \beta a_{1}$ and $z_{2}\rightarrow \beta b_{1}$, thus
\begin{equation}
 \left\{
   \begin{aligned}
   &\lambda_{+}\rightarrow \frac{1}{2}+\frac{(z_{1}-\beta a_{1})+i(z_{2}-\beta b_{1})}{2(2\beta a_{1}+4c_{1}+2i\beta b_{1})}=\frac{1}{2},\\
   &\lambda_{-}\rightarrow \frac{1}{2}-\frac{(z_{1}+\beta a_{1})+i(z_{2}+\beta b_{1})}{2(2\beta a_{1}+4c_{1}+2i\beta b_{1})}=\frac{1}{2}-\frac{\beta a_{1}+i\beta b_{1}}{2\beta a_{1}+4c_{1}+2i\beta b_{1}},\\
   \end{aligned}
   \right.\nonumber
\end{equation}
which means that for $\beta>0$, the eigenvalues of the preconditioned matrix $\mathcal{P}_{MGSSP}^{-1}\mathcal{A}$ tend to scatter near the point $(\frac{1}{2},0)$ and the point $(\frac{\beta a_{1}c_{1}+2c_{1}^{2}}{(\beta a_{1}^{2}+2c_{1})^{2}+\beta^{2}b_{1}^{2}},-\frac{\beta b_{1}c_{1}}{(\beta a_{1}^{2}+2c_{1})^{2}+\beta^{2}b_{1}^{2}})$ as $\alpha\rightarrow{0_{+}}$. Additionally, it is easily seen that the eigenvalues of $\mathcal{P}_{MGSSP}^{-1}\mathcal{A}$ tend to scatter near the points $(\frac{\alpha_{0}\beta_{0}^{2}a_{1}+2\beta_{0}^{2}(a_{1}^{2}+b_{1}^{2})+12\beta_{0} a_{1}c_{1}+(\alpha_{0}\beta_{0}+4c_{1})(4c_{1}+z_{1})+2\beta_{0}(a_{1}z_{1}+|b_{1}z_{2}|)}{2[(\alpha_{0}\beta_{0}+2\beta_{0}a_{1}+4c_{1})^{2}+4\beta_{0}^{2}b_{1}^{2}]}, \frac{(\alpha_{0}\beta_{0}+2\beta_{0} a_{1}+4c_{1})(z_{2}+\beta_{0} b_{1})-2\beta_{0} b_{1}(\beta_{0} a_{1}+4c_{1}+z_{1})}{2[(\alpha_{0}\beta_{0}+2\beta_{0}a_{1}+4c_{1})^{2}+4\beta_{0}^{2}b_{1}^{2}]})$ and $(\frac{\alpha_{0}\beta_{0}^{2}a_{1}+2\beta_{0}^{2}(a_{1}^{2}+b_{1}^{2})+12\beta_{0} a_{1}c_{1}+(\alpha_{0}\beta_{0}+4c_{1})(4c_{1}-z_{1})-2\beta_{0}(a_{1}z_{1}+|b_{1}z_{2}|)}{2[(\alpha_{0}\beta_{0}+2\beta_{0}a_{1}+4c_{1})^{2}+4\beta_{0}^{2}b_{1}^{2}]},
\frac{(\alpha_{0}\beta_{0}+2\beta_{0} a_{1}+4c_{1})(\beta_{0} b_{1}-z_{2})-2\beta_{0} b_{1}(\beta_{0} a_{1}+4c_{1}-z_{1})}{2[(\alpha_{0}\beta_{0}+2\beta_{0}a_{1}+4c_{1})^{2}+4\beta_{0}^{2}b_{1}^{2}]})$ as $\alpha\rightarrow \alpha_{0}$ and $\beta\rightarrow\beta_{0}$ $(0\leq \alpha_{0}<+\infty,0<\beta_{0}<+\infty)$. \hfill$\blacksquare$

\newtheorem{rem1}{Remark}[section]
\begin{rem1}\label{rem1}
\emph{It follows from Theorem 5.1 that
\begin{eqnarray}
Re(\lambda_{+})&=&\frac{\alpha\beta^{2}a_{1}+2\beta^{2}(a_{1}^{2}+b_{1}^{2})+12\beta a_{1}c_{1}+(\alpha\beta+4c_{1})(4c_{1}+z_{1})+2\beta(a_{1}z_{1}+|b_{1}z_{2}|)}{2[(\alpha\beta+2\beta a_{1}+4c_{1})^{2}+4\beta^{2}b_{1}^{2}]}>0,\nonumber\\
Re(\lambda_{-})&=&\frac{\alpha\beta^{2}a_{1}+2\beta^{2}(a_{1}^{2}+b_{1}^{2})+12\beta a_{1}c_{1}+(\alpha\beta+4c_{1})(4c_{1}-z_{1})-2\beta(a_{1}z_{1}+|b_{1}z_{2}|)}{2[(\alpha\beta+2\beta a_{1}+4c_{1})^{2}+4\beta^{2}b_{1}^{2}]}\nonumber\\
&\geq&\frac{8c_{1}(\beta a_{1}+2c_{1})}{2[(\alpha\beta+2\beta a_{1}+4c_{1})^{2}+4\beta^{2}b_{1}^{2}]}>0\nonumber
\end{eqnarray}
as $\alpha\geq{0}$, $\beta>0$ and $B^{T}u\neq{0}$, and if $B$ is of full column rank and $B^{T}u={0}$, then from (19), we infer that $Re(\lambda)>0$, where $(\lambda,(u^{*},v^{*})^{*})$ is an eigenpair of the preconditioned matrix $\mathcal{P}_{MGSSP}^{-1}\mathcal{A}$.
Thus all eigenvalues of $\mathcal{P}_{MGSSP}^{-1}\mathcal{A}$ have positive real parts and lie in a positive box as $B$ is of full column rank, which may result in fast convergence of Krylov subspace acceleration. Besides, from the proof of Theorem 5.1, it can be seen that when $B^{T}u=0$ and $\alpha\rightarrow{0_{+}}$, it holds that $\lambda\rightarrow \frac{1}{2}$ or $\lambda=0$; when $B^{T}u\neq0$, $\lambda\rightarrow(\frac{1}{2},0)$ as $\beta\rightarrow{0_{+}}$ for $\alpha\geq0$. This implies that the MGSSP preconditioned matrix $\mathcal{P}_{MGSSP}^{-1}\mathcal{A}$ with proper parameters $\alpha$ and $\beta$ has much denser spectrum distribution compared with the saddle point matrix $\mathcal{A}$. This means that when the MGSSP preconditioner is applied for the GMRES method, the rate of convergence (semi-convergence) can be improved considerably. This fact is further confirmed by the numerical results presented in Tables 2-4 and Tables 6-8 of Section 6. What is more, since
\begin{eqnarray}
&&(\alpha\beta+2\beta a_{1})^{2}+(\beta|b_{1}|+\sqrt{\beta^{2}b_{1}^{2}+4\alpha\beta c_{1}})^{2}\nonumber\\
&=&(\alpha\beta+2\beta a_{1})^{2}+2\beta^{2}b_{1}^{2}+4\alpha\beta c_{1}+2\beta|b_{1}|\sqrt{\beta^{2}b_{1}^{2}+4\alpha\beta c_{1}}\nonumber\\
&\leq&(\alpha\beta+2\beta a_{1})^{2}+2\beta^{2}b_{1}^{2}+4\alpha\beta c_{1}+2\beta|b_{1}|\sqrt{\beta^{2}|b_{1}|^{2}+4\alpha\beta c_{1}+\left(\frac{2\alpha c_{1}}{|b_{1}|}\right)^{2}}\nonumber\\
&=&(\alpha\beta+2\beta a_{1})^{2}+4\beta^{2}b_{1}^{2}+8\alpha\beta c_{1}\leq(\alpha\beta+2\beta a_{1}+4c_{1})^{2}+4\beta^{2}b_{1}^{2},\nonumber
\end{eqnarray}
then from (32), we have that
\begin{eqnarray}
\left|\lambda_{\pm}-\frac{1}{2}\right|^{2}\leq\frac{(\alpha\beta+2\beta a_{1})^{2}+(\beta|b_{1}|+\sqrt{\beta^{2}b_{1}^{2}+4\alpha\beta c_{1}})^{2}}{4[(\alpha\beta+2\beta a_{1}+4c_{1})^{2}+4\beta^{2}b_{1}^{2}]}\leq\frac{1}{4}
.\nonumber
\end{eqnarray}
Then all eigenvalues of $\mathcal{P}_{MGSSP}^{-1}\mathcal{A}$ are located in a circle centered at $(\frac{1}{2},0)$ with radius $\frac{1}{2}$.}
\end{rem1}

Owing to the fact the convergence of Krylov subspace methods is not only dependent on the eigenvalue distribution of the preconditioned matrix, but also on the corresponding eigenvectors of the preconditioned matrix \cite{54,55} except for the case that the preconditioned matrix is
symmetric. We next discuss the eigenvector distribution of $\mathcal{P}_{MGSSP}^{-1}\mathcal{A}$ in the following theorem.
\newtheorem{thm6}[thm1]{Theorem}
\begin{thm6}\label{thm1}
Let the MGSSP preconditioner $\mathcal{P}_{MGSSP}$ be defined as in (4). If $B$ is of full column rank and $\alpha=0$, then the preconditioned matrix $\mathcal{P}_{MGSSP}^{-1}\mathcal{A}$ has $m+i$ $(0\leq i \leq m)$ linearly independent eigenvectors, and if $B$ is of full column rank and $\alpha>0$, then the preconditioned matrix $\mathcal{P}_{MGSSP}^{-1}\mathcal{A}$ has $i$ $(0\leq i \leq m)$ linearly independent eigenvectors. If $B$ is rank deficient and $\alpha=0$, then the preconditioned matrix $\mathcal{P}_{MGSSP}^{-1}\mathcal{A}$ has $m+i+j$ $(0\leq i\leq{m},1\leq j\leq{n})$ linearly independent eigenvectors, and if $B$ is rank deficient and $\alpha>0$, then the preconditioned matrix $\mathcal{P}_{MGSSP}^{-1}\mathcal{A}$ has $i+j$ $(0\leq i\leq{m},1\leq j\leq{n})$ linearly independent eigenvectors. There are\\
1)~$m$ eigenvectors of the form $\left(
                                    \begin{array}{c}
                                      u_{l} \\
                                      0 \\
                                    \end{array}
                                  \right)
$ $(1\leq l\leq m)$ that correspond to the eigenvalue $\frac{1}{2}$ as $\alpha=0$, where $u_{l}\neq{0}$ $(1\leq l\leq m)$ are arbitrary linearly independent
vectors;\\
2)~If $B$ is of full column rank and $\alpha>0$, $i$ $(0\leq{i}\leq m)$ eigenvectors of the form $\left(
                                           \begin{array}{c}
                                             u_{l}^{1} \\
                                             \frac{(2\lambda-1)B^{T}u_{l}^{1}}{\lambda\beta} \\
                                           \end{array}
                                         \right)
$ $(1\leq l\leq i)$ that correspond to the eigenvalues $\lambda\neq\frac{1}{2}$, where $u_{l}^{1}$ $(1\leq l\leq i)$ satisfy $\lambda\beta Au_{l}^{1}=\beta\lambda^{2}(\alpha I+2A)u_{l}^{1}+(2\lambda-1)^{2}BB^{T}u_{l}^{1}$.\\
3)~If $B$ is rank deficient and $\alpha>0$, $i$ $(0\leq{i}\leq m)$ eigenvectors of the form $\left(
                                           \begin{array}{c}
                                             u_{l}^{1} \\
                                             \frac{(2\lambda-1)B^{T}u_{l}^{1}}{\lambda\beta} \\
                                           \end{array}
                                         \right)
$ $(1\leq l\leq i)$ that correspond to the eigenvalues $\lambda\neq\frac{1}{2},0$, where $u_{l}^{1}$ $(1\leq l\leq i)$ satisfy $\lambda\beta Au_{l}^{1}=\beta\lambda^{2}(\alpha I+2A)u_{l}^{1}+(2\lambda-1)^{2}BB^{T}u_{l}^{1}$; and $j$ $(1\leq j\leq{n})$ eigenvectors $\left(
                                           \begin{array}{c}
                                             0 \\
                                             v_{l}^{2} \\
                                           \end{array}
                                         \right)
$ $(1\leq{l}\leq{j})$ that correspond to the eigenvalue $0$, where $v_{l}^{2}\neq{0}$ $(1\leq{l}\leq{j})$ satisfy $Bv_{l}^{2}=0$.
\end{thm6}
\noindent {\bf Proof.}
Let $\lambda$ be an eigenvalue of the preconditioned matrix $\mathcal{P}_{MGSSP}^{-1}\mathcal{A}$ and $\left(
                                                                                                                                \begin{array}{c}
                                                                                                                                  u \\
                                                                                                                                  v \\
                                                                                                                                \end{array}
                                                                                                                              \right)
$ be the corresponding eigenvector. To investigate the eigenvector distribution, we consider Equation (24) as follows:
\begin{equation}
 \left\{
   \begin{aligned}
   &Au=\lambda(\alpha I+2A)u+(2\lambda-1) Bv,\\
   &(2\lambda-1)B^{T}u=\lambda\beta v.\\
   \end{aligned}
   \right.
\end{equation}
Firstly, we consider $B$ has full column rank. If $u=0$, then it follows from the second equation of (33) that $\lambda v=0$ and therefore $v=0$, which contradicts to the assumption that $(u^{*},v^{*})^{*}$ is an eigenvector. Hence $u\neq{0}$. If $\lambda=\frac{1}{2}$, then from (33) we can easily get $\alpha u=0$ and $v=0$. If $\alpha=0$, then Equation (33) is always true for the case of $\lambda=\frac{1}{2}$. Hence, there are $m$ linearly independent eigenvectors $\left(
                                                                                             \begin{array}{c}
                                                                                               u_{l} \\
                                                                                               0 \\
                                                                                             \end{array}
                                                                                           \right)
$ $(l = 1,2,\cdots,m)$ corresponding to the eigenvalue $\frac{1}{2}$ as $\alpha=0$, where $u_{l}$ $(l = 1,2,\cdots,m)$
are arbitrary linearly independent vectors. If $\alpha>0$, then $u=0$ and $v=0$, a contradiction.

Next, we consider the case $\lambda\neq\frac{1}{2}$. It follows from the second equation of (33) that $v=\frac{(2\lambda-1)B^{T}u}{\lambda\beta}$. Substituting this relation into the first equation of (33) results in
\begin{eqnarray}
\lambda\beta Au=\beta\lambda^{2}(\alpha I+2A)u+(2\lambda-1)^{2}BB^{T}u.
\end{eqnarray}
If there exists $u\neq{0}$ which satisfies (34), there will be $i$ $(1\leq i\leq m)$ linearly independent eigenvectors $\left(
                                                                                   \begin{array}{c}
                                                                                     u_{l}^{1} \\
                                                                                     v_{l}^{1} \\
                                                                                   \end{array}
                                                                                 \right)
$ $(1\leq{l}\leq{i})$ corresponding to the eigenvalues $\lambda\neq\frac{1}{2}$. Here, $u_{l}^{1}\neq{0}$ $(1\leq{l}\leq{i})$ satisfy $\lambda\beta Au_{l}^{1}=\beta\lambda^{2}(\alpha I+2A)u_{l}^{1}+(2\lambda-1)^{2}BB^{T}u_{l}^{1}$ and the forms of $v_{l}^{1}$ $(1\leq{l}\leq{i})$ are
\begin{eqnarray}
v_{l}^{1}=\frac{(2\lambda-1)B^{T}u_{l}^{1}}{\lambda\beta}.\nonumber
\end{eqnarray}

If $B$ is rank deficient, then $\lambda=0$ is an eigenvalue of $\mathcal{P}_{MGSSP}^{-1}\mathcal{A}$. If $\lambda=0$, then from (33), it holds that $B^{T}u=0$ and $Au=-Bv$, which lead to $B^{T}A^{-1}Bv=0$ and therefore $Bv=0$ is due to the fact that $A^{-1}$ is positive definite, then we have $Au=0$ and $u=0$. Recalling that $B$ is rank deficient, then there exists $v\neq{0}$ which satisfies $Bv=0$, hence there will be $j$ $(1\leq j\leq{n})$ linearly independent eigenvectors $\left(
                                                                                   \begin{array}{c}
                                                                                     0 \\
                                                                                     v_{l}^{2} \\
                                                                                   \end{array}
                                                                                 \right)
$ $(1\leq{l}\leq{j})$ corresponding to the eigenvalue $0$, where $v_{l}^{2}\neq{0}$ $(1\leq{l}\leq{j})$ satisfy $Bv_{l}^{2}=0$. With a quite similar strategy utilized in the case that $B$ has full column rank, we also can obtain the eigenvectors that correspond to $\lambda=\frac{1}{2}$ and $\lambda\neq{0},\frac{1}{2}$ are the same as those for the case that $B$ is of full column rank.

Now, we show that the $m+i$ eigenvectors are linearly independent when $B$ is of full column rank and $\alpha=0$. Let $c^{(1)}=[c_{1}^{(1)},c_{2}^{(1)},\cdots,c_{m}^{(1)}]$ and $c^{(2)}=[c_{1}^{(2)},c_{2}^{(2)},\cdots,c_{i}^{(2)}]$ be two vectors with $0\leq i\leq {m}$. Then, we need to show that
\begin{eqnarray}
\left(
  \begin{array}{ccc}
    u_{1} & \cdots & u_{m} \\
    0 & \cdots & 0 \\
  \end{array}
\right)\left(
         \begin{array}{c}
           c_{1}^{(1)} \\
           \vdots \\
           c_{m}^{(1)} \\
         \end{array}
       \right)
+\left(
  \begin{array}{ccc}
    u_{1}^{1} & \cdots & u_{i}^{1} \\
    v_{1}^{1} & \cdots & v_{i}^{1} \\
  \end{array}
\right)\left(
         \begin{array}{c}
           c_{1}^{(2)} \\
           \vdots \\
           c_{i}^{(2)} \\
         \end{array}
       \right)=\left(
                 \begin{array}{c}
                   0 \\
                   \vdots \\
                   0 \\
                 \end{array}
               \right)
\end{eqnarray}
holds if and only if the vectors $c^{(1)}$ and $c^{(2)}$ both are zero vectors. Recall that in (35) the first matrix arises from the case
$\lambda_{l}=\frac{1}{2}$ $(l=1,2,\cdots,m)$ in 1),
and the second matrix from the case $\lambda_{l}\neq\frac{1}{2}$ $(l=1,2,\cdots,i)$ in 2). Multiplying both sides
of (35) from left with $2\mathcal{P}_{MGSSP}^{-1}\mathcal{A}$ leads to
\begin{eqnarray}
\left(
  \begin{array}{ccc}
    u_{1} & \cdots & u_{m} \\
    0 & \cdots & 0 \\
  \end{array}
\right)\left(
         \begin{array}{c}
           c_{1}^{(1)} \\
           \vdots \\
           c_{m}^{(1)} \\
         \end{array}
       \right)
+\left(
  \begin{array}{ccc}
    u_{1}^{1} & \cdots & u_{i}^{1} \\
    v_{1}^{1} & \cdots & v_{i}^{1} \\
  \end{array}
\right)\left(
         \begin{array}{c}
           2\lambda_{1}c_{1}^{(2)} \\
           \vdots \\
           2\lambda_{i}c_{i}^{(2)} \\
         \end{array}
       \right)=\left(
                 \begin{array}{c}
                   0 \\
                   \vdots \\
                   0 \\
                 \end{array}
               \right).
\end{eqnarray}
Then, by subtracting (35) from (36), it holds that
\begin{eqnarray}
\left(
  \begin{array}{ccc}
    u_{1}^{1} & \cdots & u_{i}^{1} \\
    v_{1}^{1} & \cdots & v_{i}^{1} \\
  \end{array}
\right)\left(
         \begin{array}{c}
           (2\lambda_{1}-1)c_{1}^{(2)} \\
           \vdots \\
            (2\lambda_{i}-1)c_{i}^{(2)} \\
         \end{array}
       \right)=\left(
                 \begin{array}{c}
                   0 \\
                   \vdots \\
                   0 \\
                 \end{array}
               \right).\nonumber
\end{eqnarray}
Since the eigenvalues $\lambda_{l}\neq\frac{1}{2}$ and $\left(
                                           \begin{array}{c}
                                             u_{l}^{1} \\
                                             v_{l}^{1} \\
                                           \end{array}
                                         \right)
$ $(1\leq l\leq{i})$ are linearly independent, we infer that $c_{l}^{(2)}=0$ $(l=1,2,\cdots,i)$. Because of the linear independence of $u_{l}$ $(l=1,2,\cdots,m)$, it follows that $c_{l}^{(1)}=0$ $(l=1,2,\cdots,m)$.
Therefore, the $m+i$ eigenvectors are linearly independent.

In the sequel, we verify the $m+i+j$ eigenvectors are linearly independent when $B$ is rank deficient and $\alpha=0$. Let $c^{(1)}=[c_{1}^{(1)},c_{2}^{(1)},\cdots,c_{m}^{(1)}]$, $c^{(2)}=[c_{1}^{(2)},c_{2}^{(2)},\cdots,c_{i}^{(2)}]$ and $c^{(3)}=[c_{1}^{(3)},c_{2}^{(3)},\cdots,c_{j}^{(3)}]$ be three vectors with $0\leq i\leq {m}$ and $1\leq j\leq{n}$, and
\begin{eqnarray}
\left(
  \begin{array}{ccc}
    u_{1} & \cdots & u_{m} \\
    0 & \cdots & 0 \\
  \end{array}
\right)\left(
         \begin{array}{c}
           c_{1}^{(1)} \\
           \vdots \\
           c_{m}^{(1)} \\
         \end{array}
       \right)
+\left(
  \begin{array}{ccc}
    u_{1}^{1} & \cdots & u_{i}^{1} \\
    v_{1}^{1} & \cdots & v_{i}^{1} \\
  \end{array}
\right)\left(
         \begin{array}{c}
           c_{1}^{(2)} \\
           \vdots \\
           c_{i}^{(2)} \\
         \end{array}
       \right)+\left(
  \begin{array}{ccc}
    0 & \cdots & 0 \\
    v_{1}^{2} & \cdots & v_{j}^{2} \\
  \end{array}
\right)\left(
         \begin{array}{c}
           c_{1}^{(3)} \\
           \vdots \\
           c_{j}^{(3)} \\
         \end{array}
       \right)=\left(
                 \begin{array}{c}
                   0 \\
                   \vdots \\
                   0 \\
                 \end{array}
               \right).
\end{eqnarray}
It is necessary for us to prove that (37) holds if and only if the vectors $c^{(1)}$, $c^{(2)}$ and $c^{(3)}$ are all zero vectors, where the first matrix consists of the eigenvectors corresponding to the eigenvalue $\frac{1}{2}$ for the case 1),
and the second and the third matrices consist of those for the case 3). Premultiplying (37) with $2\mathcal{P}_{MGSSP}^{-1}\mathcal{A}$ and going through the same algebraic operations as before, we also obtain
\begin{eqnarray}
\left(
  \begin{array}{ccc}
    u_{1}^{1} & \cdots & u_{i}^{1} \\
    v_{1}^{1} & \cdots & v_{i}^{1} \\
  \end{array}
\right)\left(
         \begin{array}{c}
           (2\lambda_{1}-1)c_{1}^{(2)} \\
           \vdots \\
            (2\lambda_{i}-1)c_{i}^{(2)} \\
         \end{array}
       \right)-\left(
  \begin{array}{ccc}
    0 & \cdots & 0 \\
    v_{1}^{2} & \cdots & v_{j}^{2} \\
  \end{array}
\right)\left(
         \begin{array}{c}
           c_{1}^{(3)} \\
           \vdots \\
           c_{j}^{(3)} \\
         \end{array}
       \right)=\left(
                 \begin{array}{c}
                   0 \\
                   \vdots \\
                   0 \\
                 \end{array}
               \right).\nonumber
\end{eqnarray}
Inasmuch as $\lambda_{l}\neq\frac{1}{2}$ and $u_{l}^{1}$ $(1\leq l\leq{i})$ are linearly independent, it holds that $c_{l}^{(2)}=0$ $(l=1,2,\cdots,i)$. Then it has
\begin{eqnarray}
\left(
  \begin{array}{ccc}
    0 & \cdots & 0 \\
    v_{1}^{2} & \cdots & v_{j}^{2} \\
  \end{array}
\right)\left(
         \begin{array}{c}
           c_{1}^{(3)} \\
           \vdots \\
           c_{j}^{(3)} \\
         \end{array}
       \right)=\left(
                 \begin{array}{c}
                   0 \\
                   \vdots \\
                   0 \\
                 \end{array}
               \right).\nonumber
\end{eqnarray}
As the vectors $v_{l}^{2}$ $(l=1,2,\cdots,j)$ are also linearly
independent, we have $c_{l}^{(3)}=0$ $(l=1,2,\cdots,j)$. Thus, (37) becomes to
\begin{eqnarray}
\left(
  \begin{array}{ccc}
    u_{1} & \cdots & u_{m} \\
    0 & \cdots & 0 \\
  \end{array}
\right)\left(
         \begin{array}{c}
           c_{1}^{(1)} \\
           \vdots \\
           c_{m}^{(1)} \\
         \end{array}
       \right)=\left(
                 \begin{array}{c}
                   0 \\
                   \vdots \\
                   0 \\
                 \end{array}
               \right).\nonumber
\end{eqnarray}
Since $u_{l}$ $(l=1,2,\cdots,m)$ are linearly independent, we have $c_{l}^{(1)}=0$ $(l=1,2,\cdots,m)$.
As a result, it holds that the $m+i+j$ eigenvectors are linearly independent.

Finally, we prove the $i+j$ eigenvectors are linearly independent when $B$ is rank deficient and $\alpha>0$. Let $c^{(1)}=[c_{1}^{(1)},c_{2}^{(1)},\cdots,c_{i}^{(1)}]$ and $c^{(2)}=[c_{1}^{(2)},c_{2}^{(2)},\cdots,c_{j}^{(2)}]$ be two vectors with $0\leq i\leq {m},1\leq{j}\leq{n}$. Then, we need to show that
\begin{eqnarray}
\left(
  \begin{array}{ccc}
    u_{1}^{1} & \cdots & u_{i}^{1} \\
    v_{1}^{1} & \cdots & v_{i}^{1} \\
  \end{array}
\right)\left(
         \begin{array}{c}
           c_{1}^{(1)} \\
           \vdots \\
           c_{i}^{(1)} \\
         \end{array}
       \right)+\left(
  \begin{array}{ccc}
    0 & \cdots & 0 \\
    v_{1}^{2} & \cdots & v_{j}^{2} \\
  \end{array}
\right)\left(
         \begin{array}{c}
           c_{1}^{(2)} \\
           \vdots \\
           c_{j}^{(2)} \\
         \end{array}
       \right)=\left(
                 \begin{array}{c}
                   0 \\
                   \vdots \\
                   0 \\
                 \end{array}
               \right)\nonumber
\end{eqnarray}
holds if and only if the vectors $c^{(1)}$ and $c^{(2)}$ both are zero vectors. Since $u_{l}^{1}$ $(1\leq l\leq{i})$ are linearly independent, we infer that $c_{l}^{(1)}=0$ $(l=1,2,\cdots,i)$. Because of the linear independence of $v_{l}^{2}$ $(l=1,2,\cdots,j)$, it follows that $c_{l}^{(2)}=0$ $(l=1,2,\cdots,j)$.
Consequently, the above $i+j$ eigenvectors are linearly independent. \hfill$\blacksquare$

\section{Numerical experiments}
\begin{figure}
  \centering
  \includegraphics[width=15.5cm,height=5cm]{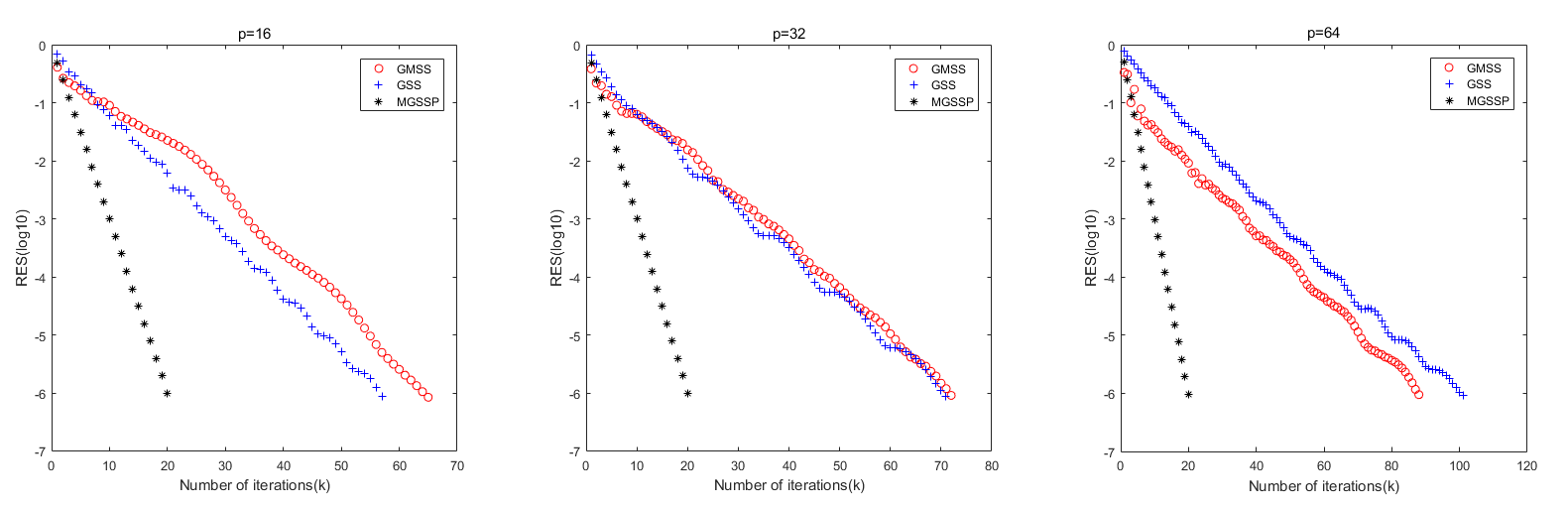}\\
  \caption{Convergence curve of algorithms with $v=0.1$ for $p=16$, $p=32$ and $p=64$, respectively.}\label{2}
\end{figure}
In this section, two numerical examples are used to verify the performance of the MGSSP iteration method and the MGSSP preconditioned GMRES method. In the meanwhile, we compare the MGSSP iteration method with the GSS and GMSS methods, and also compare the MGSSP preconditioner
with the SS, GSS, MSS and GMSS preconditioners for the GMRES method according to the number of iterations (denoted by ``IT'') and the elapsed CPU times (denoted by ``CPU''). All codes are run in MATLAB (version R2016a) and all experiments are performed on an Intel(R)
Pentium(R) CPU G3240T 2.70 GHz, 4.0GB memory and XP operating system. In our implementations, the linear
systems $(\alpha I+A+\frac{1}{\alpha}BB^{T})x=b$, $(\alpha I+A+\frac{1}{\beta}BB^{T})x=b$ and $(\alpha I+2A+\frac{4}{\beta}BB^{T})x=b$ involved in the SS, GSS and MGSSP iteration, respectively are
solved exactly by the the LU factorization. In addition, the linear
systems $(\alpha I+2H+\frac{1}{\alpha}BB^{T})x=b$ and $(\alpha I+2H+\frac{1}{\beta}BB^{T})x=b$ contained
in the MSS and the GMSS iteration are solved exactly by the Cholesky factorization.

In our numerical experiments, we choose the right-hand side vector $b$ so that the exact solution of the saddle point problem (1) is $(1,1,\cdots,1)^{T}$. All experiments are started from the initial vector $\mathbf{x}^{(0)}=(x^{(0)T},y^{(0)T})^{T}=(0,0,\cdots,0)^{T}$, terminated once the current iterate $\mathbf{x}^{(k)}$ satisfies
\begin{eqnarray}
\mathrm{RES}=\frac{\sqrt{\|f-Ax^{(k)}-By^{(k)}\|_{2}^{2}+\|g-B^{T}x^{(k)}\|_{2}^{2}}}{\sqrt{\|f\|_{2}^{2}+\|g\|_{2}^{2}}}<10^{-6},
\end{eqnarray}
and we use ``--'' to indicate that the corresponding iteration method does not satisfy the prescribed stopping criterion until $500$ iteration steps.
\label{}
\begin{Exa}\label{E1}
Consider the nonsymmetric nonsingular saddle point problem structured as (1) with the following coefficient sub-matrices \cite{43}:
\setlength\arraycolsep{4pt}
\begin{eqnarray}
&&A=\left(
    \begin{array}{cc}
      I\otimes T+T\otimes I & 0 \\
      0 & I\otimes T+T\otimes I \\
    \end{array}
  \right)\in{\mathbb{R}}^{2p^{2}\times 2p^{2}},\quad
B=\left(
    \begin{array}{c}
      I\otimes F \\
      F\otimes I \\
    \end{array}
  \right)\in{\mathbb{R}}^{2p^{2}\times p^{2}},\nonumber\\
&&T=\frac{v}{h^{2}}.\mathrm{tridiag}(-1,2,-1)+\frac{1}{2h}.\mathrm{tridiag}(-1,0,1)\in{\mathbb{R}}^{p\times{p}},\quad F=\frac{1}{h}.\mathrm{tridiag}(-1,1,0)\in{\mathbb{R}}^{p\times{p}},\nonumber
\end{eqnarray}
where $\otimes$ denotes the Kronecker product symbol and $h=\frac{1}{p+1}$ is the discretization mesh size.
\end{Exa}
\begin{figure}
  \centering
  \includegraphics[width=15.5cm,height=5cm]{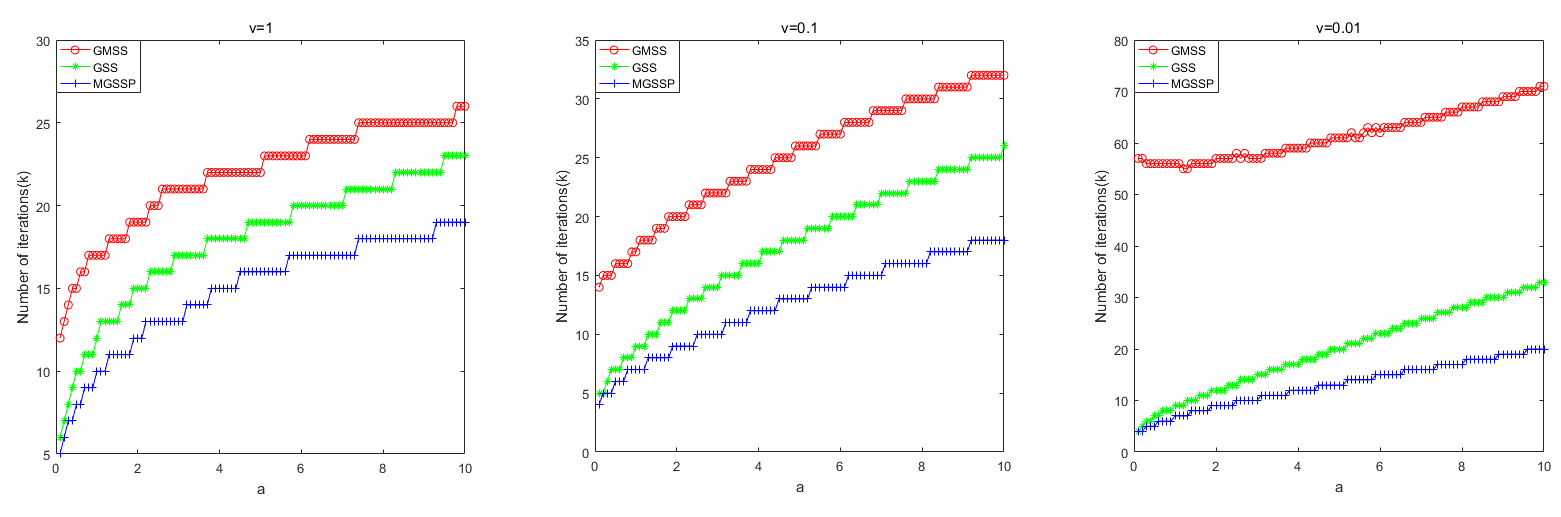}\\
  \caption{Convergence curve of algorithms with varying $\alpha=\beta$ for $p=32$.}\label{1}
\end{figure}
\begin{figure}
  \centering
  \includegraphics[width=15.5cm,height=7cm]{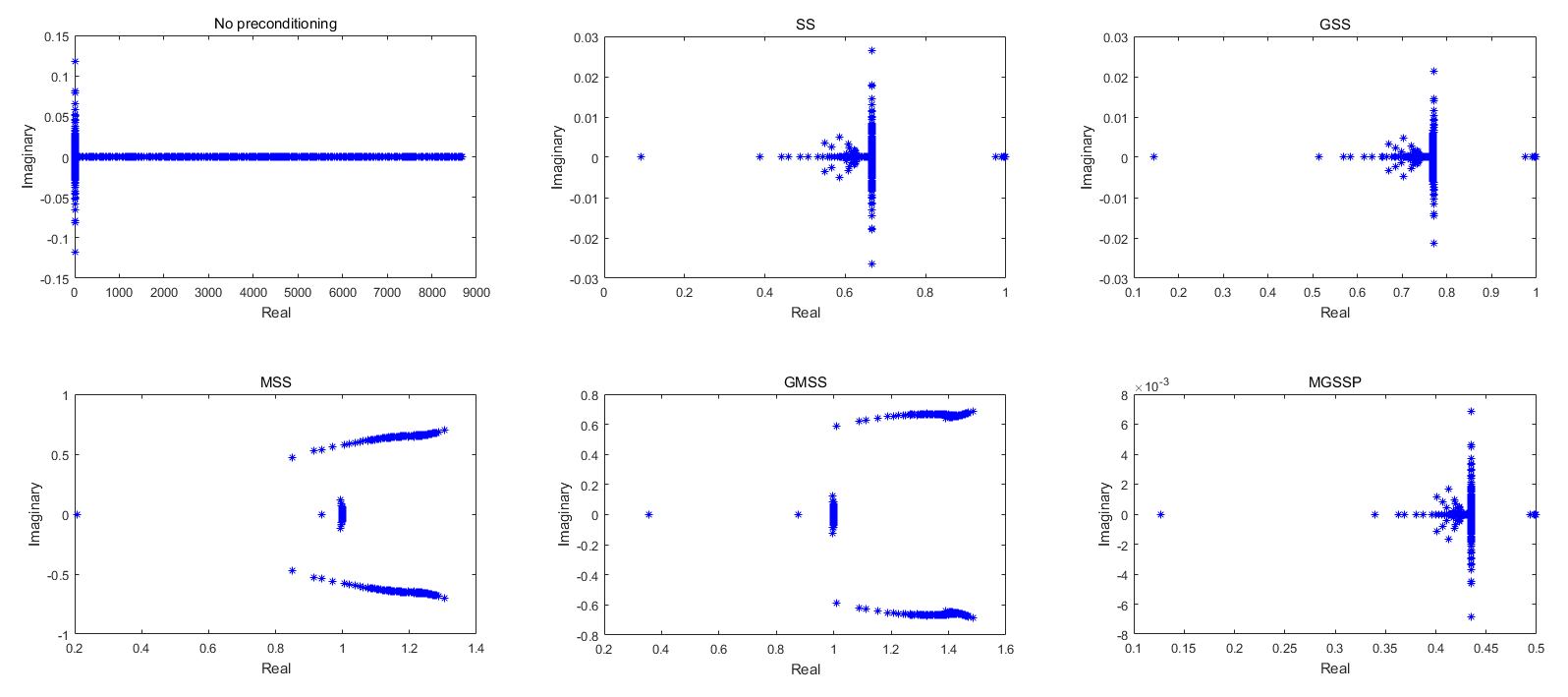}\\
  \caption{The eigenvalue distribution of the six preconditioners for $\mathcal{A}$ with $\alpha=0.6$ and $\beta=0.8$ for $p=32$ and $v=1$.}\label{1}
\end{figure}

\begin{table}[htbp]
\centering
\caption{\label{tab:test}Numerical results for the three iteration methods with $v=0.1$.}
\tabcolsep 0.3in
\begin{tabular}{lccccccccccccl}
\toprule
Method & & $p$ &  &  \\
\cline{3-5}
& &  16 & 32 & 64\\
\midrule
& $\alpha_{exp}$ & 20 & 51 & 125 \\
& $\beta_{exp}$ & 2.7 & 5 & 1.5 \\
GSS & IT & 58 & 72 & 102 \\
 &CPU & 0.2556 & 1.1583 & 21.2715\\
&RES & 8.79e-07 & 8.68e-07 & 9.80e-07 \\
& $\alpha_{exp}$ & 22 & 36 & 38 \\
& $\beta_{exp}$ & 16 & 8.3 & 5.9 \\
GMSS & IT & 66 & 73 & 89 \\
 &CPU & 0.4955 & 1.5732 & 26.4811 \\
&RES & 8.45e-07 & 9.09e-07 & 9.50e-07 \\
& $\alpha_{exp}$ & 0.2 & 0.5 & 0.2 \\
& $\beta_{exp}$ & 0.1 & 0.1 & 0.1 \\
MGSSP & IT & 21 & 21 & 21 \\
 &CPU & 0.1514 & 0.7012 & 10.1562 \\
&RES & 9.88e-07 & 9.85e-07 & 9.57e-07 \\
\bottomrule
\end{tabular}
\end{table}
In Table 1, we list the parameters involved in the tested methods which are chosen to be the experimentally found
optimal ones that minimize the total number of iteration steps for those methods, as well as the numerical results of the GSS, GMSS and MGSSP iteration methods when $v = 0.1$ with respect to different grids $16\times 16$, $32\times 32$ and $64\times 64$. Moreover, numerical results of the GMRES method and the preconditioned GMRES methods incorporated with the SS, GSS, MSS, GMSS and the MGSSP preconditioners are listed in Tables 2-4 for $v=1$, $0.1$ and $0.01$ on different uniform grids, respectively.

In order to better understand the numerical results in Table 1, convergence history of the GSS, GMSS and MGSSP iteration methods with experimental optimal parameters are depicted in Figure 1. To further confirm the effectiveness of the MGSSP preconditioned GMRES method, we plot the IT of the three preconditioned GMRES methods with parameters $\alpha=\beta$ from 0.1 to 10 with step size 0.1 in Figure 2. For more investigations, the eigenvalue distributions of the original matrix $\mathcal{A}$ and the five preconditioned matrices with $\alpha=0.6$ and $\beta=0.8$ for $v = 1$ and $p=32$ are displayed in Figure 3.

Looking into Tables 1-4 and Figures 1-3 one may make the following observations.
\begin{table}[htbp]
\centering
\caption{\label{tab:test}Numerical results for the six preconditioned GMRES methods with $v=1$, $\alpha=0.6$ and $\beta=0.8$.}
\tabcolsep 0.15in
\begin{tabular}{lcccccccccl}
\toprule
$p$ &  & $I$&  $\mathcal{P}_{SS}$ & $\mathcal{P}_{GSS}$ & $\mathcal{P}_{MSS}$ & $\mathcal{P}_{GMSS}$ & $\mathcal{P}_{MGSSP}$\\
\midrule
16& IT & 121 & 9 & 9 & 15 & 13 & 7\\
 & CPU & 0.1550 & 0.0447 & 0.1505 & 0.1838 & 0.1705 & 0.0837\\
 & RES & 7.21e-07 & 5.61e-07 & 3.29e-07 & 3.29e-07 & 6.69e-07 & 3.78e-07\\
32& IT & 264 & 10 & 9 & 15 & 14 & 7\\
 & CPU & 3.8574 & 0.5004 & 0.4703 & 0.8234 & 0.7600 & 0.3831\\
 & RES & 9.74e-07 & 2.54e-07 & 7.25e-07 & 7.63e-07 & 9.24e-07 & 9.67e-07\\
48& IT & 429 & 10 & 10 & 16 & 15 & 8\\
 & CPU & 24.7021 & 3.7594 & 3.5617 & 6.2255 & 5.7171 & 2.8951\\
 & RES & 9.95e-07 & 6.29e-07 & 1.77e-07 & 6.33e-07 & 5.21e-07 & 1.78e-07\\
64& IT & -- & 11 & 10 & 16 & 15 & 8\\
 & CPU & -- & 22.5881 & 21.3997 & 33.6562 & 31.2381 & 16.2309\\
 & RES & -- & 3.75e-07 & 2.58e-07 & 8.29e-07 & 8.93e-07 & 2.50e-07\\
\bottomrule
\end{tabular}
\end{table}

\begin{itemize}
\item{From Table 1, it can be observed that the experimental optimal parameters of the MGSSP iteration method are more stable compared with those of other two methods. Besides, the results in Table 1 imply that the MGSSP iteration method is superior to the other two methods from the point view of the IT and CPU times, and the IT of the MGSSP iteration method remains constant under the experimental optimal parameters with the increasing of the problem size.}
\end{itemize}
\begin{itemize}
\item{By comparing the results in Tables 2-4, it can be seen that without preconditioning, the GMRES method converges very slow even invalid within $500$ iteration steps for larger linear systems. All aforementioned preconditioners can largely accelerate the convergence rate of the GMRES method. The proposed MGSSP preconditioned GMRES method performs better than other five preconditioned GMRES methods as it requires less IT and CPU times. Another observation which can be pointed out is that, the convergence behavior of the MGSSP preconditioned GMRES method is not sensitive to $p$, in the sense the iterations barely change.}
\end{itemize}
\begin{table}[htbp]
\centering
\caption{\label{tab:test}Numerical results for the six preconditioned GMRES methods with $v=0.1$, $\alpha=1$ and $\beta=0.8$.}
\tabcolsep 0.15in
\begin{tabular}{lcccccccccl}
\toprule
$p$ &  & $I$&  $\mathcal{P}_{SS}$ & $\mathcal{P}_{GSS}$ & $\mathcal{P}_{MSS}$ & $\mathcal{P}_{GMSS}$ & $\mathcal{P}_{MGSSP}$\\
\midrule
16& IT & 115 & 8 & 8 & 17 & 17 & 6\\
 & CPU & 0.1326 & 0.0982 & 0.0850 & 0.2155 & 0.2371 & 0.0783\\
 & RES & 9.50e-07 & 4.54e-07 & 1.56e-07 & 5.89e-07 & 4.29e-07 & 5.96e-07\\
32& IT & 240 & 9 & 8 & 17 & 17 & 7\\
 & CPU & 3.4868 & 0.4959 & 0.4568 & 0.8974 & 0.8876 & 0.4244\\
 & RES & 9.34e-07 & 2.10e-07 & 6.49e-07 & 8.09e-07 & 4.93e-07 & 7.16e-07\\
48& IT & 367 & 9 & 9 & 18 & 17 & 7\\
 & CPU & 20.4798 & 3.3953 & 3.3951 & 6.8142 & 6.4413 & 2.7642\\
 & RES & 9.80e-07 & 4.38e-07 & 1.36e-07 & 3.97e-07 & 6.37e-07 & 1.40e-07\\
64& IT & 495 & 9 & 9 & 18 & 17 & 7\\
 & CPU & 81.8770 & 18.4334 & 18.5499 & 37.0634 & 35.4719 & 15.1190\\
 & RES & 9.73e-07 & 6.88e-07 & 2.15e-07 & 4.82e-07 & 7.42e-07 & 2.18e-07\\
\bottomrule
\end{tabular}
\end{table}
\begin{itemize}
\item{Figure 1 indicates that the three tested methods converge while the MGSSP iteration method
returns better numerical results than the GSS and the GMSS iteration methods. From Figure 2, as we expected for Example 6.1, we see that the MGSSP preconditioned
GMRES method outperforms the other two methods with the changing of $\alpha$, and show that our proposed preconditioner is more effective and practical for solving the nonsymmetric nonsingular saddle point problems, in comparison with the other preconditioners. Additionally, as seen from Figure 3, the preconditioned matrix $\mathcal{P}_{MGSSP}^{-1}\mathcal{A}$ has more clustered eigenvalues than the other ones. This means that the MGSSP preconditioner outperforms the
other five preconditioners for the GMRES method, which is congruous with the results of Table 2.}
\end{itemize}
\begin{table}[htbp]
\centering
\caption{\label{tab:test}Numerical results for the six preconditioned GMRES methods with $v=0.01$, $\alpha=1.2$ and $\beta=1.5$.}
\tabcolsep 0.15in
\begin{tabular}{lcccccccccl}
\toprule
$p$ &  & $I$&  $\mathcal{P}_{SS}$ & $\mathcal{P}_{GSS}$ & $\mathcal{P}_{MSS}$ & $\mathcal{P}_{GMSS}$ & $\mathcal{P}_{MGSSP}$\\
\midrule
16& IT & 246 & 9 & 10 & 51 & 54 & 7\\
 & CPU & 0.3743 & 0.1904 & 0.1345 & 0.5309 & 0.8654 & 0.1071\\
 & RES & 9.65e-07 & 8.26e-07 & 3.04e-07 & 9.10e-07 & 7.90e-07 & 8.61e-07\\
32& IT & 429 & 9 & 10 & 55 & 56 & 7\\
 & CPU & 7.2934 & 0.4691 & 0.5084 & 2.6743 & 3.2658 & 0.4067\\
 & RES & 9.88e-07 & 8.40e-07 & 3.19e-07 & 8.57e-07 & 8.40e-07 & 8.72e-07\\
48& IT & -- & 9 & 10 & 57 & 58 & 7\\
 & CPU & -- & 3.5469 & 3.8431 & 21.1369 & 23.6098 & 2.7461\\
 & RES & -- & 8.61e-07 & 3.21e-07 & 9.91e-07 & 9.05e-07 & 8.50e-07\\
64& IT & -- & 9 & 10 & 57 & 57 & 7\\
 & CPU & -- & 18.9807 & 20.9074 & 113.9468 & 117.1580 & 15.2981\\
 & RES & -- & 7.78e-07 & 2.58e-07 & 9.56e-07 & 9.58e-07 & 7.56e-08\\
\bottomrule
\end{tabular}
\end{table}
\begin{Exa}\label{E1}
Consider the nonsymmetric singular saddle point problem structured as (1) with the following coefficient sub-matrices \cite{60}:
\setlength\arraycolsep{4pt}
\begin{eqnarray}
A=\left(
    \begin{array}{cc}
      I\otimes T+T\otimes I & 0 \\
      0 & I\otimes T+T\otimes I \\
    \end{array}
  \right)\in{\mathbb{R}}^{2p^{2}\times{2p^{2}}},\
B=\left(
                 \begin{array}{ccc}
                   \hat{B} & b_{1} & b_{2} \\
                 \end{array}
               \right)
  \in{\mathbb{R}}^{2p^{2}\times{(p^{2}+2)}},\nonumber
\end{eqnarray}
where
\begin{eqnarray}
&&T=\frac{v}{h^{2}}.\mathrm{tridiag}(-1,2,-1)+\frac{1}{2h}.\mathrm{tridiag}(-1,0,1)\in{\mathbb{R}}^{p\times{p}},\
\hat{B}=\left(
    \begin{array}{c}
      I\otimes F \\
      F\otimes I \\
    \end{array}
  \right)\in{\mathbb{R}}^{2p^{2}\times{p^{2}}},\nonumber\\
&&b_{1}=\hat{B}\left(
               \begin{array}{c}
                 e \\
                 0 \\
               \end{array}
             \right),\
b_{2}=\hat{B}\left(
               \begin{array}{c}
                 0 \\
                 e \\
               \end{array}
             \right),\ e=(1,1,\cdots,1)\in{\mathbb{R}}^{p^{2}/2},\nonumber\\
&&F=\frac{1}{h}.\mathrm{tridiag}(-1,1,0)\in{\mathbb{R}}^{p\times{p}},\ h=\frac{1}{p+1}.\nonumber
\end{eqnarray}
Here $\otimes$ denotes the Kronecker product and $h=\frac{1}{p+1}$ is the discretization meshsize. The iterations of all tested methods are terminated once the current iterate $\mathbf{x}^{(k)}$ satisfies (38) or the maximum prescribed number of iterations $k_{max}=500$ is exceeded.
\end{Exa}
\begin{table}[htbp]
\centering
\caption{\label{tab:test}Numerical results for the three iteration methods with $v=0.1$.}
\tabcolsep 0.3in
\begin{tabular}{lccccccccccccl}
\toprule
Method & & $p$ &  &  \\
\cline{3-5}
& &  16 & 32 & 64\\
\midrule
& $\alpha_{exp}$ & 13 & 29 & 66 \\
& $\beta_{exp}$ & 39 & 53 & 60 \\
GSS & IT & 85 & 136 & 230 \\
 &CPU & 0.3211 & 2.1355 & 47.4292\\
&RES & 9.48e-07 & 9.83e-07 & 9.73e-07 \\
& $\alpha_{exp}$ & 16 & 18 & 24 \\
& $\beta_{exp}$ & 75 & 134.4 & 240 \\
GMSS & IT & 143 & 213 & 337 \\
 &CPU & 0.8641 & 4.5574 & 100.4047 \\
&RES & 9.86e-07 & 9.90e-07 & 9.98e-07 \\
& $\alpha_{exp}$ & 0.02 & 0.01 & 0.05 \\
& $\beta_{exp}$ & 0.1 & 0.05 & 0.1 \\
MGSSP & IT & 21 & 21 & 21 \\
 &CPU & 0.0990 & 0.6706 & 11.0860 \\
&RES & 9.53e-07 & 9.54e-07 & 9.54e-07 \\
\bottomrule
\end{tabular}
\end{table}

Table 5 reports the iteration counts, CPU times and relative residual
(RES) of the tested iteration methods with respect to different values of the problem size $p$ for $v=0.1$. We adopt the parameters of the tested methods to be the experimentally found optimal ones. From Table 5, we observe that although all tested methods succeed in producing approximate solutions in all cases, the MGSSP iteration method outperforms other two methods in terms of the IT and CPU times, and the advantage of the MGSSP iteration method becomes more pronounced as the system size increases.

With respect to different sizes of the coefficient matrix, we list the numerical results of the SS, GSS, MSS, GMSS and MGSSP preconditioned GMRES methods with different
values of $v$ ($v=1$, $v=0.1$ and $v=0.01$) in Tables 6-8, respectively. From Tables 6-8, we can conclude some observations as follows. Firstly, the GMRES method does not converge when $v=0.01$ and $p$ becomes large. Secondly, the five preconditioners can improve the convergence behavior of the GMRES method, but the MGSSP preconditioned GMRES method returns better numerical results than the other preconditioned GMRES methods in terms of IT and CPU time. Lastly, the MSS and GMSS preconditioned GMRES methods have worse convergence behaviors as $v$ becomes small.

The graphs of RES(log10) against number of iterations of in Table 5 for three different sizes are displayed in Figure 4. As observed in Figure 4, the MGSSP iteration method leads to much better performance than the GSS and the GMSS iteration methods. It is worthy noting that the IT of the GSS and the GMSS iteration methods increase when $p$ becomes large, but this is not true for the MGSSP iteration method.
\begin{table}[htbp]
\centering
\caption{\label{tab:test}Numerical results for the six preconditioned GMRES methods with $v=1$, $\alpha=0.6$ and $\beta=0.8$.}
\tabcolsep 0.15in
\begin{tabular}{lcccccccccl}
\toprule
$p$ &  & $I$&  $\mathcal{P}_{SS}$ & $\mathcal{P}_{GSS}$ & $\mathcal{P}_{MSS}$ & $\mathcal{P}_{GMSS}$ & $\mathcal{P}_{MGSSP}$\\
\midrule
16& IT & 145 & 9 & 8 & 15 & 13 & 6\\
 & CPU & 0.2146 & 0.0447 & 0.0968 & 0.1838 & 0.1705 & 0.0813\\
 & RES & 7.95e-07 & 5.61e-07 & 2.16e-07 & 3.29e-07 & 6.69e-07 & 8.81e-07\\
32& IT & 278 & 10 & 9 & 15 & 14 & 7\\
 & CPU & 4.1297 & 0.5004 & 0.4388 & 0.8234 & 0.7600 & 0.4251\\
 & RES & 9.79e-07 & 2.54e-07 & 1.46e-07 & 7.63e-07 & 9.24e-07 & 2.09e-07\\
48& IT & 366 & 10 & 9 & 16 & 15 & 7\\
 & CPU & 20.2558 & 3.7594 & 3.4238 & 6.2255 & 5.7171 & 2.6904\\
 & RES & 9.71e-07 & 6.29e-07 & 4.26e-07 & 6.33e-07 & 5.21e-07 & 5.66e-07\\
64& IT & 465 & 11 & 9 & 16 & 15 & 8\\
 & CPU & 76.1434 & 22.5881 & 18.1854 & 33.6562 & 31.2381 & 16.5291\\
 & RES & 9.71e-07 & 3.75e-07 & 8.37e-07 & 8.29e-07 & 8.93e-07 & 9.27e-08\\
\bottomrule
\end{tabular}
\end{table}

In order to compare effects of the GSS, GMSS, and the MGSSP preconditioned
GMRES methods in terms of the parameters $\alpha$ and $\beta$, we test these methods with $\alpha=\beta$ and plot the IT of the three preconditioned
GMRES methods with $\alpha$ from 0.1 to 10 with step size 0.1 in Figure 5. From Figure 5, we can obtain the same results as those of Figure 2.

In order to better investigate the performance of the tested preconditioned GMRES methods, Figure 6 depicts the eigenvalue distributions
of the saddle point matrix $\mathcal{A}$, the SS, GSS, MSS, GMSS and MGSSP preconditioned matrices with $v=0.1$ and $p=32$. These subfigures clearly show that the preconditioned matrices have more tightly clustered eigenvalues than the original matrix. Moreover, the eigenvalues of the MGSSP preconditioned matrix are much tighter than the other ones. These observations imply that the MGSSP preconditioned GMRES method has better numerical performance than other preconditioned GMRES methods and it can act as an efficient preconditioner for solving the nonsymmetric singular saddle point problem by the preconditioned GMRES method.

\begin{table}[htbp]
\centering
\caption{\label{tab:test}Numerical results for the six preconditioned GMRES methods with $v=0.1$, $\alpha=1.8$ and $\beta=1.5$.}
\tabcolsep 0.15in
\begin{tabular}{lcccccccccl}
\toprule
$p$ &  & $I$&  $\mathcal{P}_{SS}$ & $\mathcal{P}_{GSS}$ & $\mathcal{P}_{MSS}$ & $\mathcal{P}_{GMSS}$ & $\mathcal{P}_{MGSSP}$\\
\midrule
16& IT & 122 & 9 & 9 & 19 & 19 & 7\\
 & CPU & 0.1422 & 0.1298 & 0.1276 & 0.3788 & 0.3205 & 0.1223\\
 & RES & 8.71e-07 & 5.68e-07 & 2.67e-07 & 7.02e-07 & 5.33e-07 & 1.62e-07\\
32& IT & 237 & 10 & 9 & 19 & 19 & 7\\
 & CPU & 3.3291 & 0.6218 & 0.5728 & 1.1392 & 1.1254 & 0.4491\\
 & RES & 9.87e-07 & 2.31e-07 & 5.53e-07 & 6.67e-07 & 5.62e-07 & 2.83e-07\\
48& IT & 350 & 10 & 9 & 19 & 19 & 7\\
 & CPU & 19.3841 & 4.1896 & 3.8573 & 7.9254 & 8.3686 & 3.1636\\
 & RES & 9.99e-07 & 3.20e-07 & 7.60e-07 & 7.45e-07 & 5.74e-07 & 3.81e-07\\
64& IT & 461 & 10 & 9 & 19 & 19 & 7\\
 & CPU & 76.1620 & 23.0578 & 20.8402 & 44.0056 & 42.1010 & 17.2456\\
 & RES & 9.82e-07 & 3.91e-07 & 9.15e-07 & 8.17e-07 & 6.00e-07 & 4.58e-07\\
\bottomrule
\end{tabular}
\end{table}
\begin{figure}
  \centering
  \includegraphics[width=15.5cm,height=5cm]{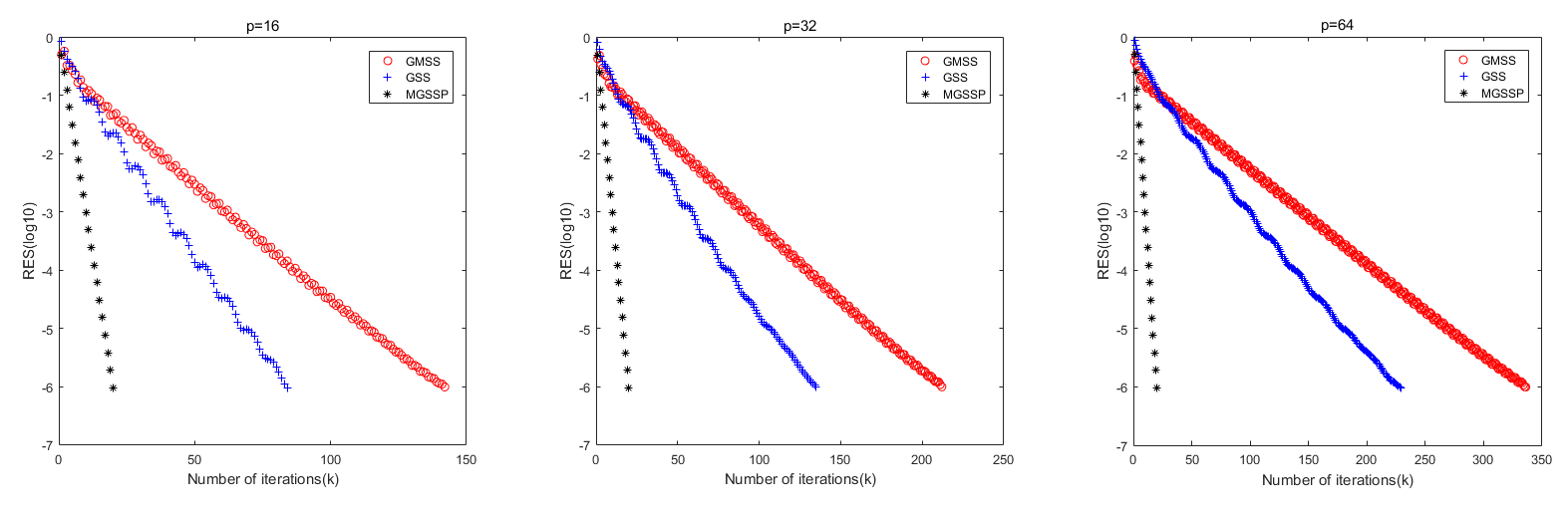}\\
  \caption{Convergence curve of algorithms with $v=0.1$ for $p=16$, $p=32$ and $p=64$, respectively.}\label{7}
\end{figure}
\begin{figure}
  \centering
  \includegraphics[width=15.5cm,height=5cm]{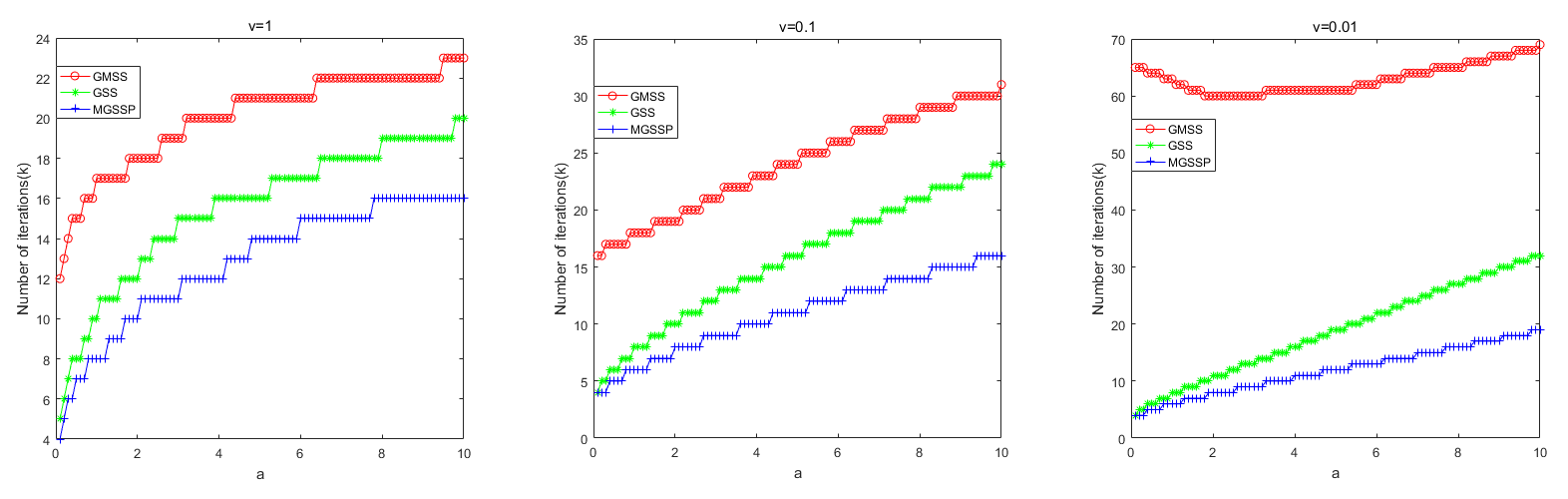}\\
  \caption{Convergence curve of algorithms with varying $\alpha=\beta$ for $p=32$.}\label{6}
\end{figure}
\begin{table}[htbp]
\centering
\caption{\label{tab:test}Numerical results for the six preconditioned GMRES methods with $v=0.01$, $\alpha=1.85$ and $\beta=1.75$.}
\tabcolsep 0.15in
\begin{tabular}{lcccccccccl}
\toprule
$p$ &  & $I$&  $\mathcal{P}_{SS}$ & $\mathcal{P}_{GSS}$ & $\mathcal{P}_{MSS}$ & $\mathcal{P}_{GMSS}$ & $\mathcal{P}_{MGSSP}$\\
\midrule
16& IT & 250 & 10 & 10 & 59 & 59 & 7\\
 & CPU & 0.4184 & 0.1320 & 0.1333 & 0.8935 & 0.9098 & 0.1191\\
 & RES & 9.42e-07 & 5.70e-07 & 4.30e-07 & 8.22e-07 & 8.03e-07 & 8.47e-07\\
32& IT & 419 & 10 & 10 & 60 & 60 & 7\\
 & CPU & 7.0626 & 0.6292 & 0.6296 & 3.3032 & 3.3277 & 0.4901\\
 & RES & 9.85e-07 & 5.31e-07 & 3.93e-07 & 9.43e-07 & 9.14e-07 & 8.11e-07\\
48& IT & -- & 10 & 10 & 60 & 60 & 7\\
 & CPU & -- & 4.1788 & 4.3580 & 23.9170 & 23.6870 & 3.0734\\
 & RES & -- & 5.75e-07 & 4.28e-07 & 9.46e-07 & 9.15e-07 & 8.72e-07\\
64& IT & -- & 10 & 10 & 60 & 60 & 7\\
 & CPU & -- & 23.6671 & 24.2199 & 130.4922 & 130.1184 & 16.9088\\
 & RES & -- & 6.43e-07 & 4.83e-07 & 9.40e-07 & 9.07e-07 & 9.34e-07\\
\bottomrule
\end{tabular}
\end{table}
\section{Conclusions}
For nonsymmetric saddle point problems, by combining the GSS and MSSP of a matrix, we establish a modified generalized shift-splitting (MGSSP) iteration method and the corresponding preconditioner called the MGSSP preconditioner in this paper. The unconditional convergence and semi-convergence of the MGSSP iteration method for solving nonsingular and singular saddle point problems, respectively are discussed. Moreover, eigenproperties of the preconditioned matrix are described. Numerical results given in Section 6 illustrate that the efficiency of the MGSSP iteration method and the MGSSP preconditioner for saddle point problems with nonsymmetric positive definite (1,1) parts, and confirm that they outperform some existing ones.

\begin{figure}
  \centering
  \includegraphics[width=15.5cm,height=7cm]{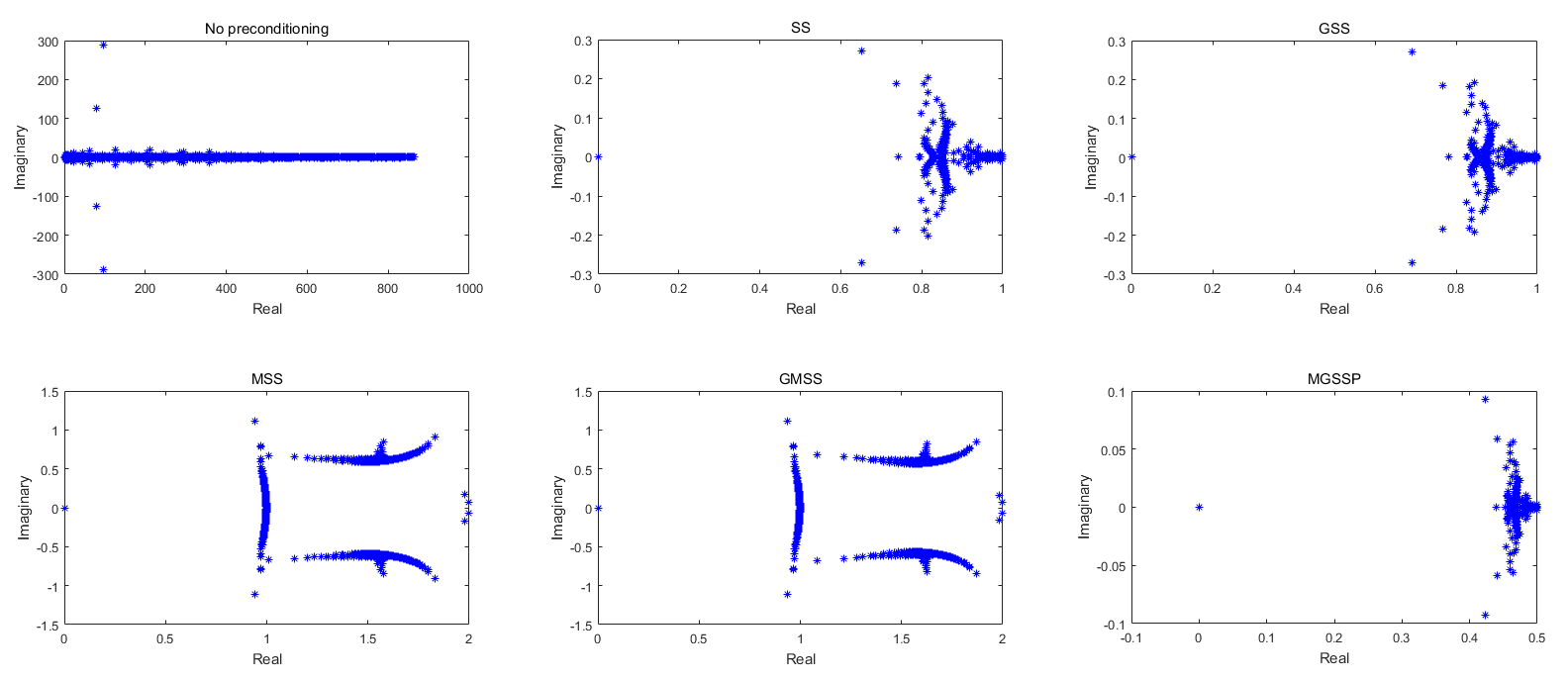}\\
  \caption{The eigenvalue distribution of the six preconditioners for $\mathcal{A}$ with $\alpha=1.8$ and $\beta=1.5$ for $p=32$ and $v=0.1$.}\label{8}
\end{figure}
We should point out that the MGSSP preconditioner may not have the optimality property, i.e., the iteration counts depend on the parameters $\alpha$ and $\beta$ (see Figures 2 and 5). Besides, admittedly, the choices of the optimal parameters of the MGSSP iteration method and the MGSSP preconditioned GMRES method is a challenging problem that deserves further study. For most iterative methods, this work is very complicated. Nevertheless, by adopting certain approximation strategies, there have been practically useful formula for obtaining nearly optimal iteration parameters; see \cite{19,46,47}. To further investigations, we would like to study how to further improve the MGSSP preconditioner and choose the optimal parameters for the MGSSP iteration method.



\end{document}